\documentclass[11pt]{article}
\usepackage[active]{srcltx}
\usepackage{amssymb,amsmath}

\newtheorem{proposition}{Proposition}[section]
\newtheorem{theorem}{Theorem}[section]
\newtheorem{lemma}{Lemma}[section]
\newtheorem{corollary}{Corollary}[section]
\newtheorem{definition}{Definition}[section]
\newtheorem{remark}{Remark}[section]

\numberwithin{equation}{section}

\def\H{\mathbb{H}}

\def\R{\mathbb{R}}

\def\L{\mathcal{L}}
\def\eps{\varepsilon}

\begin{document}

\begin{center}{\Large\sc
Sharp comparison and maximum principles via horizontal normal mapping in
the Heisenberg group
}\\
\vspace{0.2cm} {\large Zolt\'an M. Balogh\footnote{Z. M. Balogh was
supported by the Swiss National Science Foundation, and the FP7 EU
Commission Project CG-DICE.}, Andrea Calogero, and Alexandru
Krist\'aly \footnote{A. Krist\'aly  was supported by a CNCS-UEFISCDI
grant
no. PN-II-ID-PCE-2011-3-0241, and  J\'anos Bolyai Research
Scholarship.
}}
\end{center}

\begin{abstract}
\noindent  In this paper we solve a problem raised by Guti\'errez
and Montanari about comparison principles  for $H-$convex functions
on subdomains of Heisenberg groups. Our approach is based on the
notion of the sub-Riemannian horizontal normal mapping and uses
degree theory for set-valued maps. The statement of the comparison
principle combined with a Harnack inequality is applied to prove the
Aleksandrov-type maximum principle, describing the correct boundary
behavior of continuous $H-$convex functions vanishing at the
boundary of horizontally bounded subdomains of Heisenberg groups.
This result answers a question by Garofalo and Tournier. The
sharpness of our results are illustrated by examples.
\end{abstract}


\medskip


\tableofcontents

\section{Introduction}
\subsection{Motivation}
It is well known that convex functions defined on subdomains of
$\R^{n}$ are locally Lipschitz continuous and almost everywhere
twice differentiable. Moreover, the celebrated maximum principle due
to Aleksandrov provides a global regularity result for convex
functions that are continuous on the closure and are vanishing on
the boundary of the domain.  More precisely, if $\Omega\subset
\mathbb R^n$ is a bounded open and convex domain, and $u\in
C(\overline \Omega)$ is convex with $u=0$ on $\partial \Omega$, then
\begin{equation}\label{elso-Alex}
    |u(\xi_0)|^n\leq C_n {\rm dist}(\xi_0,\partial \Omega){\rm diam}(\Omega)^{n-1} \L^{n}(\partial
u(\Omega)),\ \forall  \xi_0\in\Omega,
\end{equation}
where $C_n>0$ is a constant depending only on the dimension $n$.  In
the above expression the notation $\L^{n}(\partial u(\Omega))$
stands for the measure of the range of the so-called normal mapping
of $u$. To define this concept we need first the subdifferential
$\partial u(\xi_0)$ of $u$ at $\xi_0$, given by
$$\partial u(\xi_0)=\left\{p\in \mathbb R^n:u(\xi)\geq u(\xi_0)+p\cdot (\xi-\xi_0),\
\forall \xi\in \Omega\right\},$$
 where
$'\cdot'$ is the usual inner product in $\mathbb R^{n}$.
The range of the normal mapping of $u$ is defined by $$\partial
u(\Omega)=\bigcup_{\xi\in\Omega}\partial u(\xi).$$ A convenient way
to deduce the Aleksandrov estimate (\ref{elso-Alex})
 is to compare the ranges of normal
mappings of the convex function $u$ and the cone function
$v:\overline\Omega\to \mathbb R$ with base on $\partial \Omega$ and
vertex $(\xi_0,u(\xi_0))$ (see e.g. Guti\'errez \cite[Theorem
1.4.2]{Gu2001}).

It is  well-known, that for any convex function $u\in C^2(\Omega)$,
\begin{equation}\label{monge-ampere-measure}
\L^{n}(\partial u(\Omega))=\int_\Omega {\rm det} [{\rm
Hess}(u)(x)]dx,
\end{equation}
which implies by \eqref{elso-Alex} the estimate:
\begin{equation}\label{maso-Alex}
    |u(\xi_0)|^n\leq C_n {\rm dist}(\xi_0,\partial \Omega){\rm diam}(\Omega)^{n-1}\int_\Omega {\rm det} [{\rm Hess}(u)(x)]dx,\ \forall  \xi_0\in\Omega.
\end{equation}
In recent years, the notion of convexity has been considered in the
setting of  Heisenberg groups by Lu, Manfredi and Stroffolini
\cite{GuMaSt2004}, and in more general Carnot groups by Danielli,
Garofalo and Nhieu \cite{DaGaNh2003} and also Juutinen, Lu, Manfredi
and Stroffolini \cite{JuGuMaSt2007}. The main idea behind this
approach is to develop a concept of convexity that is adapted to the
sub-Riemannian, or Carnot-Carath\'eodory geometry of the Carnot
groups. In this way convexity is assumed only along {\it
trajectories of left-invariant horizontal vector-fields} which are
in the first layer of the Lie algebra of the group and generate the
sub-Riemannian metric. This notion is called by many authors as {\it
$H-$convexity}. This approach makes sense also in case of more
general Carnot-Carath\'eodory spaces even in the absence of  a
groups structure, see Bardi and Dragoni \cite{Bardi-Dragoni}.

Various results on local regularity properties  such as local
Lipschitz continuity or second differentiability a.e. in terms of
the horizontal vector-fields have been already proven in this
context. We refer to the paper of
 Balogh and Rickly \cite{Balogh-Rickly} for the proof of  the local Lipschitz continuity of $H-$convex functions on the Heisenberg group and Rickly \cite{Ri2006} for Carnot groups. 
 
 It was pointed out to us by one of the referees, that  the generalization of Aleksandrov's second order differentiability theorem of $H-$convex functions to the case of Carnot groups is a rather delicate issue. Magnani \cite{Magnani} proved second horizontal differentiability a.e. in the general Carnot setting of a $H-$convex function $u$, but only under the assumption that all entries of the symmetrized horizontal Hessian $u_{i,j}$ as well as the horizontal commutators $[X_{i}, X_{j}]u$ are Radon measures. The first condition was proved by Danielli, Garofalo and Nhieu in \cite{DaGaNh2003}. The second condition is more difficult, it was proven by Danielli, Garofalo, Nhieu and Tournier in \cite{DaGaNhTo2004} for the case of Carnot groups of step 2.  The property that $[X_{i}, X_{j}]u$ are Radon measures is still open for general Carnot groups.

 In this paper we will be concerned with first order regularity properties of $H-$convex functions on the Heisenberg group. We note first, that the behavior of $H-$convex functions in non-horizontal directions can still be pretty wild. Indeed,
 examples of $H-$convex functions are constructed by Balogh and Rickly in \cite{Balogh-Rickly} which coincide with the Weierstrass function on a thick Cantor set of vertical lines.
 This fact indicates the intricate nature of $H-$convex functions as well as possible differences with respect to their Euclidean counterpart.
In particular, the validity of an Aleksandrov-type estimate, similar
to \eqref{elso-Alex} becomes questionable.

The main goal of this paper is to prove {\it global regularity
results} akin to \eqref{elso-Alex} in the setting of general
Heisenberg groups $\H^n$. This problem has been first considered by
Guti\'errez and Montanari \cite{GuMo-CommPDE} in the setting of the
first Heisenberg group $\H^{1}$ and
 by Garofalo and Tournier
\cite{GaTo2005} for the second Heisenberg group $\H^2$ and the Engel
group. In these papers, the methods of Trudinger and Wang
\cite{TrudWang1, TrudWang2, TrudWang3} have been applied to obtain
comparison estimates for integrals involving Hessians and related
expressions in second order derivatives. Trudinger and Zhang
\cite{TrudZhang} obtained recently  a generalization of these
results for integrals of $k$-th order Hessian measures  of
$k$-convex functions defined on $\H^{n}$. Such comparison estimates
can be used to deduce weaker versions of Aleksandrov-type maximum
principle \eqref{maso-Alex}. For instance, in \cite{GuMo-CommPDE} it
is shown that if $u: B_{H}\to \R$ is a $C^{2}-$smooth, $H-$convex
function defined on the unit Kor\'anyi-Cygan ball in the first
Heisenberg group $\mathbb H^1$ which vanishes on the boundary, then
 \begin{equation}\label{GuMo}|u(\xi_0)|^2\leq c_1(\xi_0) \int_{B_H} \left({\rm det} [{\rm
 Hess}_H(u)(\xi)]^*+12(Tu(\xi))^2\right)d\xi, \
 \forall \xi_0\in B_H,\end{equation}
  where $[{\rm
Hess}_H(u)(\xi)]^*$ denotes the symmetrized horizontal Hessian and
$Tu$ is the vertical derivative of $u$.

The main drawback of the estimate (\ref{GuMo}) is that the
expression $c_1(\xi_0)>0$ in front of the integral behaves like
${\rm dist}_H(\xi_0,\partial B_H)^{-\alpha}$ for some $\alpha>0$,
which is far to be optimal taking into account that $u=0$ on
$\partial B_H$. A similar result was obtained also in
\cite{GaTo2005}, where Garofalo and Tournier \cite[p.
2013]{GaTo2005}  formulated the question about existence of  a
suitable pointwise estimate that behaves like a positive power of
the distance to the boundary.

\subsection{Statements of main results }
The primary goal of our paper is to provide a positive answer to the
above question by proving an Aleksandrov-type estimate in the spirit
of (\ref{elso-Alex}). More precisely, we shall prove the estimate
\begin{equation}\label{alex-intro}
    |u(\xi_0)|^{2n}\leq C_n{\rm dist}_H(\xi_0,\partial \Omega)\
{\rm diam}_{HS}(\Omega)^{2n-1}\L_{HS}^{2n}(\partial_Hu(\Omega)),\
\forall \xi_0\in \Omega,
\end{equation}
where $\Omega\subset \mathbb H^n$ is any open horizontally bounded
and convex domain, $u:\overline \Omega\to \mathbb R$ is a continuous
$H-$convex function which vanishes at the boundary $\partial
\Omega$, and $C_n>0$ depends only on $n$. (The concept of horizontal
boundedness will be introduced in the sequel.)

 In the above estimate ${\rm dist}_H$ stands for the sub-Riemannian distance of the Heisenberg
group. The quantities ${\rm diam}_{HS}(\Omega)$ and
$\L_{HS}^{2n}(\partial_Hu(\Omega))$ denote the {\it horizontal
slicing diameter} of  the horizontally bounded set  $\Omega$, resp.
the {\it horizontal slicing measure} of the set
$\partial_Hu(\Omega)$. These notions are introduced in Definition
\ref{horiz-slice-diam} as the appropriate substitutes for their
Euclidean counterparts  ${\rm diam}(\Omega)$ and $\L^{n}(\partial
u(\Omega))$, respectively.

We recall that  $\partial_H u$ is  the  {\it horizontal normal
mapping} of $u$ introduced by Danielli, Garofalo and Nhieu
\cite{DaGaNh2003} and studied by Calogero and Pini \cite{CaPi2011}.
The  concept of horizontal normal mapping turns out to be the right
analogue to the normal mapping in the Euclidean space which made the
estimate \eqref{alex-intro} possible. Roughly speaking, the
horizontal normal mapping $\partial_{H}u$
 includes all  subdifferentials of
$u$ taken in the directions of the left-invariant  horizontal
directions on the Heisenberg group.

Until now, there was a major obstacle in applying the method of
normal mapping due to the lack of good comparison principles for
$H-$convex functions. Our first result overcomes this obstacle, and
at the same time answers a question of Calogero and Pini
\cite{CaPi2011} and Guti\'errez and Montanari \cite{GuMo-CommPDE}:

\begin{theorem}[Comparison
principle for the horizontal normal
mapping]\label{comparison-theorem-new} Let $\Omega \subset \H^{n}$
be an open, horizontally bounded and $H-$convex set, and
$u,v:\Omega\to\R$ be $H-$convex functions.   Let ${\Omega_0}\subset
\mathbb H^n$ be open such that $\overline {\Omega_0}\subset\Omega$
and assume that $u< v$ in ${\Omega_0}$ and $u=v$ on $\partial
{\Omega_0}$. Then
$$\partial_{H}v({\Omega_0}) \subset
\partial_{H}u({\Omega_0}).$$
\end{theorem}

In fact, Theorem \ref{comparison-theorem-new} is a consequence of a
more general comparison result, see Theorem \ref{n-comparison-uj},
where the novelty of our approach is shown by the application of a
degree theoretical argument for upper semicontinuous set-valued
maps, developed by  Hu and Papageorgiou \cite{Hu-Papa}. Due to the
$H-$convexity of the functions $u$ and $v$, the upper semicontinuous
set-valued maps $\partial_H u$ and $\partial_H v$ show certain
monotonicity properties, allowing to relate the set-valued degree of
these maps via a suitable homotopy flow. A similar comparison
principle to the previous one can be stated by requiring
 $u \leq v$ in $\Omega_0$ but adding the  strict $H-$ convexity of $v$, see
Theorem \ref{cor-strict}.

We emphasize that the $H-$convexity of the functions $u$ and $v$ is
indispensable in order to obtain comparison principles. Indeed, in
the absence of convexity we construct an example for which the
comparison principle fails on the first Heisenberg group $\mathbb
H^1$, see Section \ref{sect-5}.

Using Theorem \ref{comparison-theorem-new} we can prove the
following:

\begin{theorem}[Horizontal comparison principle]\label{n-comparison intro}
\noindent Let $\Omega \subset \H^{n}$ be an open, bounded and
$H-$convex set, and
 $u,v: \overline{\Omega} \to \R$ be continuous  $H-$convex functions.
 If for every Borel set $E\subset\Omega$ we have
$$
\L^{2n}(\partial_H v(E))\le \L^{2n}(\partial_H u(E)),
$$
then
$$
\min_{\xi\in\overline{\Omega}}(v(\xi)-u(\xi))=\min_{\xi\in\partial\Omega}(v(\xi)-u(\xi)).
$$
\end{theorem}

A consequence of Theorem \ref{n-comparison intro} is the fact that
the horizontal normal mapping characterizes uniquely the $H-$convex
functions with prescribed boundary values.

\begin{corollary} \label{corollary
n-comparison bis convex} \noindent Let $\Omega \subset \H^{n}$ be an
open, bounded and $H-$convex set, and let
  $u,v:\overline\Omega\to \mathbb R$ be  continuous $H-$convex functions.
 If for every Borel set $E\subset\Omega$ we have
$$
\L^{2n}(\partial_H u(E))=\L^{2n}(\partial_H v(E))
$$
and $u=v$ in $\partial\Omega,$ then $u=v$ in $\Omega.$
\end{corollary}


The main result of the paper is the following maximum principle.

\begin{theorem}[Aleksandrov-type maximum principle] \label{Alexandrov-n}
 Let $\Omega \subset \H^{n}$ be an open, horizontally  bounded and
convex set. If $u:\overline{\Omega} \to \R$ is a continuous
$H-$convex function which verifies $u=0$ on $\partial \Omega$, then
\begin{equation} \label{alexandrov bis}
 |u(\xi_0)|^{2n}\leq C_n{\rm dist}_H(\xi_0,\partial \Omega)\
{\rm
diam}_{HS}(\Omega)^{2n-1}\L_{HS}^{2n}(\partial_Hu(\Omega)),\qquad
\forall \xi_0\in \Omega,
\end{equation}
where $C_n>0$ depends only on $n$.
\end{theorem}

The proof of Theorem \ref{Alexandrov-n} is a puzzle which is
assembled by several pieces: basic comparison principle,  maximum
principle on horizontal planes, horizontal normal mapping of cone
functions, Harnack-type inequality, and quantitative description of
the twirling effect of horizontal planes. Some of the pieces in this puzzle are 
readily available in the current literature: in particular the Harnack-type inequality 
for $H-$convex functions has been proven by Gutierrez and Montanari in 
\cite{GuMo-CommPDE}, in the same paper the authors apply this result to 
obtain estimates on the boundary behavior of $H-$convex functions.

Theorem \ref{Alexandrov-n} is {\it sharp} which is shown as follows:
for a given $\eps\in (0,1)$ we construct an open, bounded and convex
set $\Omega\subset \H^{1}$ and a continuous $H-$convex function
$u:\overline{\Omega} \to \R$  which verifies $u=0$ on $\partial
\Omega$ and $u<0$ in $\Omega$ such that
$\L_{HS}^{2}(\partial_Hu(\Omega))<\infty$, and
\begin{equation}\label{limsup-hater-0}
    \sup_{\xi\in \Omega} \frac{|u(\xi)|^{2}}{{\rm dist}_H(\xi,\partial
    \Omega)^{1+\eps}}=+\infty.
\end{equation}

Some comments concerning further perspectives are in order. Since
the arguments in the proof of the comparison principles (see
Theorems \ref{comparison-theorem-new} and \ref{cor-strict}) are
topological, it is clear that such results can be also extended to
general Carnot groups. However, in this general setting certain
technical difficulties will arise in the proof of the
Aleksandrov-type maximum principle, e.g. the construction of
specific cone functions; these issues will be considered in the
forthcoming paper \cite{forth}. Furthermore, we expect that the
approach presented in this paper can be successfully applied to
establish interior $\Gamma^{1+\alpha}-,$ or $W^{2,p}-$regularity of
$H-$convex functions in the spirit of Caffarelli \cite{Caf-1, Caf-2}
and Guti\'errez \cite{Gu2001}. In the setting of Carnot groups a
first step in this direction has been done by Capogna and Maldonado
\cite{Capogna-Maldonado}.

The paper is organized as follows. In Section 2 we fix notations and
recall preliminary results on $H-$convex functions in the Heisenberg
group. Section 3 is devoted to comparison principles; in particular
we prove Theorems \ref{comparison-theorem-new} and \ref{n-comparison
intro}. In Section 4 we give the proof of our main result Theorem
\ref{Alexandrov-n}. Section \ref{sect-5} is devoted to the
discussions related to sharpness of our results. First, we provide
an example showing that comparison principles do not hold in the
absence of the convexity assumption, see \S\ref{sect-51}. Then, the
above example (see (\ref{limsup-hater-0})) is presented in detail,
showing the sharpness of the Aleksandrov-type estimate, see
\S\ref{sect-52}. We also discuss the relationship bet\-ween the
horizontal Monge-Amp\`ere operator and the horizontal normal
mapping, see \S\ref{sect-53}.  To make the  paper self-contained we 
add an  Appendix 
containing two parts. In the first part we recall those results of Hu and Papageorgiu
\cite{Hu-Papa} on
the degree theory for set-valued maps from which we need in our proof in
Section 3. In the second part of the Appendix we give a detailed proof of the
quantitative Harnack inequality following Gutierrez and Montanari \cite{GuMo-CommPDE} 
 that we use in Section 4.\\

\noindent {\it Acknowledgment.}  A. Calogero and A. Krist\'aly are
grateful to the Mathematisches Institute, Universit\"at Bern for the
warm hospitality where this work has been initiated. We thank the 
referees for the careful reading of the manuscript  and for their comments.

\section{Preliminaries}\label{sect-2}

The Heisenberg group $\mathbb H^n$ is the simplest Carnot group of
step 2 which serves as prototype of Carnot groups. For a
comprehensive introduction to analysis on Carnot groups we refer to
\cite{BLU}. Here we recall just the necessary notation and
background results used in the sequel. The  Lie algebra  $\frak{h}$
of $\mathbb H^n$ admits a stratification  $ \frak{h}=V_1\oplus V_2$
with $V_1=\texttt{\rm span}\{X_{i},\, Y_i;\ 1\le i\le n\}$ being the
first layer, and  $V_2=\texttt{\rm span}\{T\}$ being the second
layer which is one-dimensional. We assume $[X_i,Y_i]=-4T$  and the
rest of commutators of basis vectors all vanish. The exponential map
$\exp:\frak{h}\to \mathbb H^n$ is defined in  the usual way. By
these commutator rules we obtain, using the Baker-Campbell-Hausdorff
formula, that $\H^n=\mathbb C^n\times \mathbb R$ is endowed with the
non-commutative group law given by
\begin{equation}\label{Heisenberg-law}
    (z,t)\circ (z',t')=(z+z',t+t'+2 {\rm Im}  \langle z,  {z'}\rangle),
\end{equation}
where $z=(z_1,...,z_n)\in \mathbb C^n$, $t\in \mathbb R$, and
$\langle z,  {z'}\rangle=\sum_{j=1}^n z_j\overline{z_j'}$ is the
Hermitian inner product. Denoting by $z_j=x_j+iy_j$, then
$(x_1,...,x_n,y_1,...,y_n,t)$ form a real coordinate system for
$\mathbb H^n$. Transporting the basis vectors of $V_1$  from the
origin to an arbitrary point of the group by left- translations, we
obtain a system of left-invariant  vector fields written as first
order differential operators as follows
\begin{equation} \label{vector fields} \left.
  \begin{array}{lll}
     X_j=\partial_{x_j}+2y_j \partial_t,\qquad j=1,...,n;\\
     Y_j=\partial_{y_j}-2x_j\partial_t,\qquad j=1,...,n.
  \end{array}
\right.
\end{equation}
These vector fields are called by an abuse of language {\it
horizontal}.  The  {\it horizontal plane} in $\xi_0\in \mathbb H^n$
is given by  $H_{\xi_{0}}=\xi_0\circ \exp(V_1\times \{0\})$.
  It is easy to check that
for  $\xi_0=(z_0,t_0)=(x_0,y_0,t_0)\in \mathbb H^n$ the equation of
the horizontal plane is given by
$$H_{\xi_0}=\{(z,t)\in \mathbb H^n:t=t_0+2{\rm Im}\langle z_0,
{z}\rangle\}=\{(x,y,t)\in \mathbb H^n:t=t_0+2(x\cdot y_0-x_0\cdot
y)\}.$$ The sub-Riemannian, or Carnot-Carath\'eodory metric on
$\H^n$ is defined in terms of the above vector fields. Instead of
the Carnot-Carath\'eodory metric,  in this paper we shall work with
the bi-Lipschitz equivalent {\it Kor\'anyi-Cygan metric} that is
more suitable for concrete calculations and is defined explicitly as
follows.

Let $N(z,t)=(|z|^4+t^2)^\frac{1}{4}$ be the gauge norm on $\mathbb
H^n$. It is an interesting exercise to check that the expression
$$d_H((z,t),(z',t'))=N((z',t')^{-1}\circ (z,t)),$$
satisfies the triangle inequality defining a metric on $\H^{n}$ (see
\cite{Cygan}). This metric is the so-called Kor\'anyi-Cygan metric
which is by left- translation and dilation invariance bi-Lipschitz
equivalent to the Carnot-Carath\'eodory metric. Here, the
non-isotropic Heisenberg dilations $\delta_{\lambda}: \H^{n}\to
\H^{n}$ for $\lambda > 0$ are defined by $\delta_{\lambda}(z,t) =
(\lambda z, \lambda^{2}t)$. If $A\subset \mathbb H^n$ and $\xi\in
\mathbb H^n$, then ${\rm dist}_H(\xi,A)=\inf_{\zeta\in
A}d_H(\xi,\zeta)$.  The Kor\'anyi-Cygan  ball of center
$(z_0,t_0)\in \mathbb H^n$ and radius $r>0$ is given by
$B_H((z_0,t_0),r)=\{(z,t)\in \mathbb H^n:d_H((z,t),(z_0,t_0))<r\}.$

Let $\Omega\subset \mathbb H^n$ be an open set.   The main idea of
the analysis on the Heisenberg group is that general regularity
properties of functions defined on the Heisenberg group should be
expressed only in terms of horizontal vector fields \eqref{vector
fields}. In particular, the appropriate gradient notion for a
function is the so-called {\it horizontal gradient}, which is
defined as the $2n$--vector $ \nabla_Hu(\xi)=
\left(X_1u(\xi),...,X_nu(\xi), Y_1u(\xi),...,Y_nu(\xi)\right)$ for a
function $u\in\Gamma^1(\Omega)$. 
Here, the class $\Gamma^k(\Omega)$ is the Folland--Stein space of
functions having continuous derivatives up to order $k$ with respect
to the vector fields $X_i$ and $Y_i,$ $i\in \{1,...,n\}$. For
general non-smooth functions $u:\Omega\to \mathbb R$ one defines the
{\it horizontal subdifferential} $\partial _H u(\xi_0)$ of $u$ at
$\xi_0\in \Omega$
 given by
$$\partial _H u(\xi_0)=\left\{p\in \mathbb R^{2n}:u(\xi)\geq u(\xi_0)+p\cdot
({\texttt{\rm Pr}_1}(\xi)-{\texttt{\rm Pr}_1}(\xi_0)), \ \forall
\xi\in \Omega\cap H_{\xi_0}\right\},$$  where ${\texttt{\rm
Pr}_1}:\H^n\to\R^{2n}$ is the projection defined by ${\texttt{\rm
Pr}_1}(\xi)={\texttt{\rm Pr}_1}(x,y,t)=(x,y)$. (The same notation
$'\cdot'$ will be used for the inner products in $\mathbb R^n$ and
$\mathbb R^{2n}$.)
 It is easy to see that if $u\in \Gamma^1(\Omega)$ and
$\partial_Hu(\xi)\neq \emptyset$, then $\partial_Hu(\xi)=\{\nabla_H
u(\xi)\}$.

The range of the {\it horizontal normal mapping} of the function $u$
is defined by
$$\partial_H u(\Omega)=\bigcup_{\xi\in\Omega}\partial_H u(\xi).$$

  A function $u:\Omega\to \mathbb R$ is called $H-${\it
subdifferentiable on} $\Omega$ if $\partial_H u(\xi)\neq \emptyset$
for every $\xi\in \Omega.$ Let $\mathcal S_H(\Omega)$ be the set of
all $H-$subdifferentiable functions on $\Omega$, and $\mathcal
S_H^0(\Omega)$  be set of all continuous $H-$subdifferentiable
functions on $\Omega$.

The main objects of study in this paper are $H-$ convex functions.
There are several equivalent ways to define the concept of $H-$
convexity. The most intuitive property is to require  the convexity
of the restriction of the function on the trajectories of left
invariant vector fields spanned by \eqref{vector fields}. Another
definition using the group operation is as follows. A set $\tilde
\Omega\subset \mathbb H^n$ is called $H-${\it convex} if for every
$\xi_1,\xi_2\in \tilde \Omega$ with $\xi_1\in H_{\xi_2}$ and
$\lambda\in [0,1]$, we have $\xi_1\circ \delta_\lambda
(\xi_1^{-1}\circ \xi_2)\in \tilde \Omega$. It is clear that if
$\tilde \Omega$ is convex (i.e. it is convex in $\R^{2n+1}-$sense),
then it is also $H-$convex. If $\tilde \Omega$ is $H-$convex, a
function $u:\tilde \Omega\to \mathbb R$ is called $H-${\it convex}
if for every $\xi_1,\xi_2\in \tilde \Omega$ with $\xi_1\in
H_{\xi_2}$ and $\lambda\in [0,1]$, we have
\begin{equation}\label{convex-def}
    u(\xi_1\circ \delta_\lambda (\xi_1^{-1}\circ \xi_2))\leq
(1-\lambda)u(\xi_1)+\lambda u(\xi_2).
\end{equation}
 If the  strict inequality holds in (\ref{convex-def}) for every
$\xi_1\neq \xi_2,$ $ \xi_1\in H_{\xi_2}$  then $u$ is called {\it
strictly} $H-${\it convex}. We denote by  $\mathcal C_H(\tilde
\Omega)$ the set of all $H-$convex functions on $\tilde \Omega$.

We will now present some basic properties of $H-$convex functions
which will be used through the paper. First, for various equivalent
characterizations of $H-$convex functions and their regularity
properties we refer to
\cite{Balogh-Rickly,CaPi2011,DaGaNh2003,CaPi2012} which can be
summarized as follows:
\begin{theorem}\label{teorema mappa normale}
Let  $\Omega\subset \mathbb H^n$ be an open set. If $u:\Omega\to
\mathbb R$ is a function, then $\partial_H u(\xi)$ is a convex and
compact set of $\mathbb R^{2n}$ for every $\xi\in \Omega$. If
$\Omega$ is $H-$convex, then $\mathcal S_H(\Omega)=\mathcal
S_H^0(\Omega)= \mathcal C_H(\Omega)$.
\end{theorem}

 Now, we are dealing with the regularity of
the set-valued map $\xi \mapsto
\partial_{H}u(\xi)$. Let us recall that if $X$ and $Y$ are metric spaces,  a set-valued
map $F:X\to 2^Y\setminus\{\emptyset\}$ with compact values is {\rm
upper semicontinuous} at $x\in X$ if for every $\varepsilon>0$ there
exists $\delta>0$ such that for every $x'\in B_X(x,\delta)$ one has
$F(x')\subset B_Y(F(x),\varepsilon)$. $F$ is upper semicontinuous on
$Z\subset X$ if it is upper semicontinuous at every point $x\in Z.$
Here,  $B_X(x,\delta)$ and $B_Y(y,\delta)$ denote the balls of radii
$\delta$ and center $x$ and $y$, respectively, in $X$ and $Y$.

\begin{proposition}\label{usco}
Let $\Omega\subset \mathbb H^n$ be an open set. If $u\in \mathcal
S_H^0(\Omega)$  then $\partial_H u:\Omega\to 2^{\mathbb R^{2n}}$ is
{\it upper semicontinuous} on $\Omega$. Moreover, for every compact
set $K\subset \Omega$, the set $\partial_H u(K)$ is compact.

\end{proposition}

{\it Proof.} 
Let $\xi_0\in \Omega$ be fixed and assume that $\partial_H u$ is not
upper semicontinuous at $\xi_0.$ On account of the upper
semicontinuity  and Theorem \ref{teorema mappa normale}  this
implies the existence of a sequence $\{\xi_k\}\subset \Omega$ such
that $\xi_k\to\xi_0$ and $p_k\in
\partial_Hu(\xi_k)$ with $p_k\to p_0$ and $p_0\notin
\partial_Hu(\xi_0)$. Note that $p_k\in
\partial_Hu(\xi_k)$ is equivalent to $$u(\zeta)-u(\xi_k)\geq p_k\cdot (\texttt{\rm Pr}_1(\zeta)-\texttt{\rm Pr}_1(\xi_k)),\qquad \forall \zeta\in \Omega\cap H_{\xi_k}.$$
Let $\zeta\in \Omega\cap H_{\xi_0}$ be a given point and take a
sequence $\zeta_{k}\in \Omega\cap H_{\xi_k}$ with $\zeta_{k} \to
\zeta$. Then
$$u(\zeta_{k})-u(\xi_k)\geq p_k\cdot (\texttt{\rm Pr}_1(\zeta_{k})-\texttt{\rm Pr}_1(\xi_k)).$$

Since $u$ is continuous, taking the limit in the above inequality,
we have
$$u(\zeta)-u(\xi_0)\geq p_0\cdot (\texttt{\rm Pr}_1(\zeta)-\texttt{\rm Pr}_1(\xi_0)).$$
Since $ \zeta\in \Omega\cap H_{\xi_0}$ was arbitrary we obtain that
$p_0\in \partial_Hu(\xi_0)$, a contradiction. The second statement
follows (see \cite[Proposition 1.1.3]{Aubin-Cellina}) from the upper
semicontinuity of the map $\partial_H u$.\hfill $\square$\\

In the statement of our main result Theorem \ref{Alexandrov-n} the
notions of horizontal slicing  diameter ${\rm diam}_{HS}(\Omega)$
and horizontal slicing measure have been used. Roughly speaking,
${\rm diam}_{HS}(\Omega)$ stands for the supremum of diameters of
horizontal slices of $\Omega$ and
$\L_{HS}^{2n}(\partial_Hu(\Omega))$ is the supremum of measures for
the ranges of horizontal slices under the normal map. The precise
definition is as follows:
\begin{definition} \label{horiz-slice-diam}
An open  set $\Omega \subset \mathbb H^n$ is called {\rm
horizontally bounded} if
\begin{equation} \label{horiz-slice}
{\rm diam}_{HS}(\Omega)= \sup\{  {\rm diam}_{H}( \Omega\cap
H_{\xi})):\xi \in \Omega \}<+\infty.
\end{equation}
The quantity ${\rm diam}_{HS}(\Omega)$ is called the {\rm horizontal
slicing diameter} of $\Omega$. For a function $u:\Omega \to \R$ we
define the {\it horizontal slicing measure} by
$$\L_{HS}^{2n}(\partial_Hu(\Omega))=\sup_{\xi\in \Omega} \L^{2n}(\partial_{H}u(\Omega\cap H_{\xi})).$$
\end{definition}

It is clear that the quantity $ {\rm diam}_{HS}(\Omega)$ is smaller
than the Heisenberg diameter of $\Omega$ and that
$\L_{HS}^{2n}(\partial_Hu(\Omega))\leq
\L^{2n}(\partial_Hu(\Omega))$.  Theorem \ref{Alexandrov-n} implies
therefore the weaker estimate
\begin{equation} \label{alexandrov weak}
 |u(\xi_0)|^{2n}\leq C_n{\rm dist}_H(\xi_0,\partial \Omega)\
{\rm diam}_{H}(\Omega)^{2n-1}\L^{2n}(\partial_Hu(\Omega)),\qquad
\forall \xi_0\in \Omega.
\end{equation}
Notice also that $ {\rm diam}_{HS}(\Omega)$ could be finite for
certain unbounded domains $\Omega \subset \mathbb H^n$, e.g., a
cylinder around the vertical axis. Moreover, one can easily check
that we have a natural scaling invariance property of Theorem
\ref{Alexandrov-n} with respect to Heisenberg dilations
$\delta_{\lambda}$; see Remark \ref{scaling invariance}.

We conclude this section by stating some properties of $H-$convex
functions which are vanishing at the boundary.

\begin{proposition}\label{prop-non-positive-H-convex} Let $\Omega \subset \H^{n}$ be an open,  horizontally bounded and $H-$convex set. If
$u:\overline{\Omega} \to \R$ is an $H-$convex function which
verifies $u=0$ on $\partial \Omega$, then $u\leq 0$. Moreover, if
$\Omega$ is $($Euclidean$)$ convex, either $u\equiv 0$ on
$\overline{\Omega}$, or $u<0$ in $\Omega$.
\end{proposition}

{\it Proof.} Let $\xi_0\in \Omega$ be fixed. Let us consider
arbitrarily a point $\xi\in \partial\Omega\cap H_{\xi_0}.$ Since
$\Omega$ is horizontally bounded and $H-$convex, there exists a
unique point $\xi'\in (\partial\Omega\cap H_{\xi_0}\cap
H_{\xi})\setminus\{\xi\}$ such that $\xi_0=\xi\circ
\delta_\lambda(\xi^{-1}\circ\xi')$ for some
$\lambda\in(0,1)$. 
The $H-$convexity of $u:\overline{\Omega} \to \R$ implies that
$$
u(\xi_0)\le(1-\lambda)u(\xi)+\lambda u(\xi')=0,
$$
which proves that $u\le 0$ in $\Omega.$ For the proof of the second
statement we show that any two points can be connected by a certain
chain of balls where we can apply a Harnack-type inequality; we
postpone this construction to the Appendix (see Subsection
\ref{appendix-2}).
 \hfill
$\square$\\

\section{Comparison principles in Heisenberg groups}\label{sect-3}

Let us recall that in order to prove the Aleksandrov-type estimate
(\ref{elso-Alex}) in the Euclidean case, the following result is
applied (see Guti\'errez
\cite[Lemma 1.4.1]{Gu2001}):\\

\begin{lemma}[Comparison lemma in Euclidean case]\label{lemma Gu} Let
$\Omega\subset \mathbb R^n$ be an open and bounded set.  If $u,v\in
C(\overline \Omega)$ with $u=v$ on $\partial \Omega$ and $u\leq v$
in $\Omega$, then $\partial v(\Omega)\subset
\partial u(\Omega)$.
\end{lemma}

It is natural to ask whether a similar property holds in the setting of Heisenberg groups:\\
{\bf Question:} {\it Let $\Omega\subset\mathbb H^n$ be an open and
bounded set, $u,v\in C(\overline \Omega)$ with $u=v$ on $\partial
\Omega$ and $u\leq v$ in $\Omega$. Does the inclusion $\partial_H
v(\Omega)\subset
\partial_H u(\Omega)$ hold?}\\

 The answer to this question is {\it negative} in general; we postpone our counterexample to Section
 \ref{sect-5}. However, we can give  a {\it positive} answer to the
 Question formulated above, under the assumption of $H-${\it convexity}.
\subsection{Comparison lemma for the horizontal normal mapping 
}\label{sect-32}

The main result of this section is a Heisenberg version of Lemma
\ref{lemma Gu}.  While in the Euclidean case the proof of this
comparison principle is rather trivial, the geometric structure of
the Heisenberg group $\H^n$ causes serious difficulties in the proof
of such a comparison result. Various authors including  Guti\'errez
and Montanari expressed their doubts about this method and used
another approach to obtain Aleksandrov-type estimates
\cite{GuMo-CommPDE}. Here we overcome the difficulties by using
degree-theoretical arguments of set valued maps \cite{Hu-Papa}; the
results needed in the proof are collected in the Appendix. Our first
result is the following:

\begin{theorem}[Comparison lemma for horizontal normal mapping]\label{n-comparison-uj}
\noindent Let $\Omega_0$ and $\Omega \subset \H^{n}$ be open,
horizontally  boun\-ded  sets such that $\Omega$ is $H-$convex,
$\overline{\Omega_0}\subset \Omega$ and $u,v:\Omega \to \mathbb R$
are $H-$convex functions.
 Let $\xi_0 \in {\Omega_0} $ be fixed such that $u(\xi_0) \leq v(\xi_0)$ and $u \geq v$ on
$\partial {\Omega_0} \cap H_{\xi_0}$. If  $p_0\in\partial_H
v(\xi_0)$ satisfies
\begin{eqnarray}\label{4-1-relacio}
&&v(\xi)>v(\xi_0)+p_0\cdot(\texttt{\rm Pr}_{1}(\xi)-\texttt{\rm
Pr}_{1}(\xi_{0})),\qquad \forall\xi\in\partial{\Omega_0}\cap
H_{\xi_0},
\end{eqnarray}
then $p_0\in \partial_H u({\Omega_0}\cap H_{\xi_0}).$
\end{theorem}

 {\it Proof.} The proof is divided into four steps.

 {\bf Step 1.} We consider the restriction of the standard projection $\texttt{\rm Pr}_{1}$ to a horizontal plane: more precisely, consider $\texttt{\rm
Pr}_1: H_{\xi_0} \to \R^{2n}$ which gives a linear isomorphism
between the horizontal plane $H_{\xi_0}$ and $\R^{2n}$. Accordingly,
we introduce the following notations, $\tilde{\xi_0}:= \texttt{\rm
Pr}_1(\xi_0)$, $\tilde{\xi}:= \texttt{\rm Pr}_1(\xi)$,
$\widetilde{\partial_H v}:=
\partial_H v \circ \texttt{\rm Pr}_1^{-1}:\texttt{\rm Pr}_1(\overline {\Omega_0} \cap H_{\xi_0}) \to
2^{\mathbb R^{2n}}$ and $\widetilde{\partial_H u}:=
\partial_H u \circ \texttt{\rm Pr}_1^{-1}:\texttt{\rm Pr}_1(\overline {\Omega_0} \cap H_{\xi_0}) \to
2^{\mathbb R^{2n}}$. In these notations the condition
\eqref{4-1-relacio} reads as
\begin{eqnarray}\label{4-1-terelacio}
&&v(\xi)>v(\xi_0)+p_0\cdot(\tilde{\xi}-\tilde{\xi_{0}}),\qquad
\forall\xi\in\partial{\Omega_0}\cap H_{\xi_0}.
\end{eqnarray}
By Proposition \ref{usco} and Theorem \ref{teorema mappa normale},
the set-valued maps $\widetilde{\partial_H u}$ and
$\widetilde{\partial_H v}$ are upper semicontinuous on the compact
set $\texttt{\rm Pr}_1(\overline {\Omega_0} \cap H_{\xi_0})$ with
compact and convex values.

{\bf Step 2.} Let $p_0\in \partial_Hv(\xi_0).$ We prove that
\begin{equation}\label{elso-egy-sv-degree}
\texttt{\rm deg}_{SV}\left(\widetilde{\partial_H
v}(\cdot)-p_0,\texttt{\rm Pr}_1( {\Omega_0} \cap
H_{\xi_0}),0\right)=1,
\end{equation}
where $\texttt{\rm deg}_{SV}$ denotes the degree function for
set-valued maps, see Theorem \ref{Hu-Papa-theorem} from the
Appendix.

To verify (\ref{elso-egy-sv-degree}),  we first claim that
\begin{equation}\label{monotonitas-uj}
(p^v-p_0)\cdot(\tilde \xi-\tilde \xi_0)>0,\ \forall \xi\in \partial
{\Omega_0}\cap H_{\xi_0},\qquad \forall p^v\in \widetilde{\partial_H
v}(\tilde \xi).
\end{equation}
Let us fix $\xi\in
\partial {\Omega_0}\cap H_{\xi_0}$ and  $p^v\in \widetilde{\partial_H
v}(\tilde \xi)$. Since $\xi\in \Omega$ and $v$ is $H-$convex on
$\Omega$, one has that
$$ v(\zeta)-v(\xi)\geq p^v\cdot (\tilde \zeta-\tilde \xi),\qquad \forall
\zeta\in  \Omega\cap H_{\xi}.$$ In particular, choosing
$\zeta=\xi_0\in {\Omega_0}\cap H_{\xi_0}$ in the latter inequality,
we obtain that
\begin{equation}\label{kell-kesobb-v}
    v(\xi_0)-v(\xi)\geq p^v\cdot (\tilde \xi_0-\tilde \xi).
\end{equation}
Combining this inequality with (\ref{4-1-terelacio}), it yields
precisely  relation (\ref{monotonitas-uj}).

Now, we consider the parametric set-valued map $\mathcal
F_\lambda:\texttt{\rm Pr}_1(\overline {\Omega_0}\cap H_{\xi_0})\to
2^{\mathbb R^{2n}}$, $\lambda\in [0,1]$, defined by
$$\mathcal F_\lambda(\tilde \xi)=(1-\lambda)(\tilde \xi-\tilde \xi_0)+\lambda (\widetilde{\partial_H v}(\tilde \xi)-p_0).$$
It follows from Proposition \ref{usco} and Theorem \ref{teorema
mappa normale} that the following properties hold:
\begin{itemize}
  \item $\overline{\{{\cup \mathcal F_\lambda(\tilde \xi):(\lambda,\tilde \xi)\in [0,1]\times \texttt{\rm Pr}_1(\overline {\Omega_0}\cap H_{\xi_0})\}}}$ is compact in $\mathbb R^{2n}$;
  \item for every $(\lambda,\tilde \xi)\in [0,1]\times \texttt{\rm Pr}_1(\overline {\Omega_0}\cap H_{\xi_0})$, the set $\mathcal F_\lambda(\tilde \xi)$ is compact and convex in
  $\mathbb R^{2n}$;
  \item $(\lambda,\tilde \xi)\mapsto \mathcal F_\lambda(\tilde \xi)$ is upper semicontinuous from $[0,1]\times \texttt{\rm Pr}_1(\overline {\Omega_0}\cap H_{\xi_0})$ into $2^{\mathbb R^{2n}}\setminus \{\emptyset\}$.
\end{itemize}
According to Definition \ref{P-class-homo} from the Appendix,
$\mathcal F_\lambda$ is of homotopy of class (P).

We now claim that for the constant curve $\gamma:[0,1]\to \mathbb
R^{2n}$, $\gamma(\lambda)=0$, we have  $\gamma(\lambda)\notin
\mathcal F_\lambda(\texttt{\rm Pr}_1(\partial {\Omega_0}\cap
H_{\xi_0}))$ for every $\lambda\in [0,1]$. By contrary, we assume
that there exists $\lambda_0\in [0,1]$ and $\xi\in
\partial {\Omega_0}\cap H_{\xi_0}$ such that $0\in \mathcal F_{\lambda_0}(\tilde
\xi)$, i.e.,  $$0\in (1-\lambda_0)(\tilde \xi-\tilde
\xi_0)+\lambda_0 (\widetilde{\partial_H v}(\tilde \xi)-p_0).$$ In
particular, there exists $p^v\in \widetilde{\partial_H v}(\tilde
\xi)$ such that $0=(1-\lambda_0)(\tilde \xi-\tilde \xi_0)+\lambda_0
(p^v-p_0)$. Multiplying the latter relation by $(\tilde \xi-\tilde
\xi_0)\neq 0$, on account of (\ref{monotonitas-uj}) we obtain the
contradiction
$$0=(1-\lambda_0)|\tilde \xi-\tilde \xi_0|^2+\lambda_0
(p^v-p_0)\cdot(\tilde \xi-\tilde \xi_0)>0.$$ Therefore, by the
homotopy invariance (see Theorem \ref{Hu-Papa-theorem} from the
Appendix), we have that $\lambda\mapsto \texttt{\rm
deg}_{SV}(\mathcal F_\lambda,\texttt{\rm Pr}_1({\Omega_0}\cap
H_{\xi_0}),0)$ is constant. In particular, by exploiting the basic
properties of the set-valued and Brouwer degrees (see Appendix), it
yields that
\begin{eqnarray*}
 && \texttt{\rm deg}_{SV}\left(\widetilde{\partial_H v}-p_0,\texttt{\rm Pr}_1( {\Omega_0} \cap
H_{\xi_0}),0\right) = \texttt{\rm deg}_{SV}(\mathcal F_1,\texttt{\rm
Pr}_1( {\Omega_0} \cap
H_{\xi_0}),0) =\\
&&\qquad= \texttt{\rm deg}_{SV}(\mathcal F_0,\texttt{\rm Pr}_1(
{\Omega_0} \cap H_{\xi_0}),0) =\texttt{\rm deg}_{SV}(Id-\tilde
\xi_0,\texttt{\rm Pr}_1( {\Omega_0}
\cap H_{\xi_0}),0)=\\
&&\qquad=\texttt{\rm deg}_{B}(Id-\tilde \xi_0,\texttt{\rm Pr}_1(
{\Omega_0} \cap H_{\xi_0}),0)
   =  \texttt{\rm deg}_B(Id,\texttt{\rm Pr}_1( {\Omega_0} \cap H_{\xi_0}),\tilde
\xi_0)
   =1,
\end{eqnarray*}
 which shows (\ref{elso-egy-sv-degree}).

{\bf Step 3.} We prove that
$$\texttt{\rm deg}_{SV}\left(\widetilde{\partial_H u}-p_0,\texttt{\rm Pr}_1( {\Omega_0} \cap
H_{\xi_0}),0\right)=1.$$ First of all, a similar reason as in
 (\ref{kell-kesobb-v}) shows that
\begin{equation}\label{kell-kesobb-u}
    u(\xi_0)-u(\xi)\geq p^u\cdot (\tilde \xi_0-\tilde \xi),\ \forall
\xi\in \partial {\Omega_0}\cap H_{\xi_0},\qquad  \forall p^u\in
\widetilde{\partial_H u}(\tilde \xi).
\end{equation}

We introduce  the parametric set-valued map $\mathcal
G_\lambda:\texttt{\rm Pr}_1(\overline {\Omega_0}\cap H_{\xi_0})\to
2^{\mathbb R^{2n}}$, $\lambda\in [0,1]$, defined by
$$\mathcal G_\lambda(\tilde \xi)={(1-\lambda)(\widetilde{\partial_H v}(\tilde \xi)-p_0)+\lambda (\widetilde{\partial_H u}(\tilde \xi)-p_0)}.$$
We observe, again  from Proposition \ref{usco} and Theorem
\ref{teorema mappa normale} that
\begin{itemize}
  \item $\overline{\{{\cup \mathcal G_\lambda(\tilde \xi):(\lambda,\tilde \xi)\in [0,1]\times \texttt{\rm Pr}_1(\overline {\Omega_0}\cap H_{\xi_0})\}}}$ is compact in $\mathbb R^{2n}$;
  \item for every $(\lambda,\tilde \xi)\in [0,1]\times \texttt{\rm Pr}_1(\overline {\Omega_0}\cap H_{\xi_0})$,  $\mathcal G_\lambda(\tilde \xi)$ is compact and convex in
  $\mathbb R^{2n}$ (as the sum of two compact and convex sets);
  \item $(\lambda,\tilde \xi)\mapsto \mathcal G_\lambda(\tilde \xi)$ is upper semicontinuous from $[0,1]\times \texttt{\rm Pr}_1(\overline {\Omega_0}\cap H_{\xi_0})$ into $2^{\mathbb R^{2n}}\setminus \{\emptyset\}$.
\end{itemize}
Therefore, $\mathcal G_\lambda$ is a homotopy of class (P).

We prove that
\begin{equation}\label{gamma-nulla}
    0\notin \mathcal G_\lambda(\texttt{\rm Pr}_1(\partial {\Omega_0}\cap
H_{\xi_0})),\qquad \forall \lambda\in [0,1].
\end{equation}
Assume the contrary, i.e., there exists $\lambda_0\in [0,1]$ and
$\xi\in
\partial {\Omega_0}\cap H_{\xi_0}$ such that $0\in\mathcal  G_{\lambda_0}(\tilde
\xi)$. It follows that
\begin{equation}\label{0-homotop}
    0=(1-\lambda_0)(p^v-p_0)+\lambda_0 (p^u-p_0)
\end{equation}
for some  $p^u\in \widetilde{\partial_H u}(\tilde \xi)$ and $\
p^v\in \widetilde{\partial_H v}(\tilde \xi)$.
 Combining (\ref{kell-kesobb-v}),
(\ref{kell-kesobb-u}) and (\ref{0-homotop}) respectively,  we obtain
that
$$(1-\lambda_0)v(\xi_0)+\lambda_0 u(\xi_0)-[(1-\lambda_0)v( \xi)+\lambda_0 u( \xi)]\geq p_0\cdot(\tilde \xi_0-\tilde \xi).$$
On the other hand, by adding the latter inequality to
(\ref{4-1-terelacio}) applied for $\tilde \xi$, it yields
$$\lambda_0(-v(\xi_0)+u(\xi_0))+\lambda_0(v( \xi)- u( \xi))>0.$$
Note that $u\geq v$ on $\partial {\Omega_0}\cap H_{\xi_0}$; thus it
follows that
$$\lambda_0(u(\xi_0)-v(\xi_0))>\lambda_0(u( \xi)- v( \xi))\geq 0.$$
Clearly, $\lambda_0\neq 0$; thus, it yields that $u(\xi_0)>v(\xi_0)$
which contradicts the assumption that $v(\xi_0)\geq u(\xi_0)$.
Therefore, (\ref{gamma-nulla}) holds true.

Again,  by the homotopy invariance (see Theorem
\ref{Hu-Papa-theorem} from the Appendix), we have that
$\lambda\mapsto \texttt{\rm deg}_{SV}(\mathcal G_\lambda,\texttt{\rm
Pr}_1({\Omega_0}\cap H_{\xi_0}),0)$ is constant, i.e., according to
Step 2,
$$\texttt{\rm deg}_{SV}\left(\widetilde{\partial_H u}-p_0,\texttt{\rm Pr}_1( {\Omega_0} \cap
H_{\xi_0}),0\right)=\texttt{\rm deg}_{SV}\left(\widetilde{\partial_H
v}-p_0,\texttt{\rm Pr}_1( {\Omega_0} \cap H_{\xi_0}),0\right)=1,$$
which concludes the proof of Step 3.

{\bf Step 4.}  By Step 3 and the definition of $\texttt{\rm
deg}_{SV}$, for small $\varepsilon>0$, one has that
\begin{equation}\label{deg-utolso}
    \texttt{\rm deg}_B(f_\varepsilon^u-p_0,\texttt{\rm Pr}_1( {\Omega_0} \cap H_{\xi_0}),0)=1,
\end{equation}
where $f^u_\varepsilon:\texttt{\rm Pr}_1( \overline{\Omega_0} \cap
H_{\xi_0})\to \mathbb R^{2n}$ is a continuous approximate selector
of the upper semicontinuous set-valued map $\widetilde{\partial_H
u}$
 such that
\begin{equation}\label{eps-approx}
    f_\varepsilon^u(\tilde \xi)\in \widetilde{\partial_H u}\left(B_{\mathbb
R^{2n}}(\tilde \xi,\varepsilon)\cap \texttt{\rm Pr}_1( \overline
{\Omega_0} \cap H_{\xi_0})\right)+B_{\mathbb
R^{2n}}(0,\varepsilon),\ \ \forall \tilde \xi\in \texttt{\rm
Pr}_1(\overline {\Omega_0} \cap H_{\xi_0}),
\end{equation}
see Proposition \ref{prop-cellina} from the Appendix.  Let
$\varepsilon=\frac{1}{k}$ and let $\phi_k^u:=f_{1/k}^u$, $k\in
\mathbb N.$ First of all, from (\ref{deg-utolso}) and the properties
of the Brouwer degree $d_B$ (see Theorem \ref{degree-Browder} from
the Appendix), we have that for every $k\in \mathbb N$ there exists
$\tilde \xi_k\in \texttt{\rm Pr}_1( {\Omega_0} \cap H_{\xi_0})$ such
that $p_0=\phi_k^u(\tilde \xi_k)$. Up to a subsequence, we may
assume that $\tilde \xi_k\to \tilde \xi\in \texttt{\rm Pr}_1(
\overline{\Omega_0} \cap H_{\xi_0})$. On the other hand, by relation
(\ref{eps-approx}), we have that
$$p_0=\phi_k^u(\tilde \xi_k)\in \widetilde{\partial_H u}\left(B_{\mathbb R^{2n}}\left(\tilde \xi_k,\frac{1}{k}\right)\cap \texttt{\rm Pr}_1( \overline{\Omega_0}
\cap H_{\xi_0})\right)+B_{\mathbb
R^{2n}}\left(0,\frac{1}{k}\right),$$ i.e., there exists $\tilde
\zeta_k\in B_{\mathbb R^{2n}}(\tilde \xi_k,\frac{1}{k})\cap
\texttt{\rm Pr}_1(\overline {\Omega_0} \cap H_{\xi_0})$ and $p_k\in
B_{\mathbb R^{2n}}(0,\frac{1}{k})$ such that $p_0\in
\widetilde{\partial_H u}(\tilde\zeta_k)+p_k.$ Clearly,
$\tilde\zeta_k\to \tilde \xi$ as $k\to \infty$. In the following, we
shall show that $p_0\in \widetilde{\partial_H u}(\tilde \xi).$

We assume by contradiction, that $p_0\notin \widetilde{\partial_H
u}(\tilde \xi).$ Since $\widetilde{\partial_H u}(\tilde \xi)$ is
compact, it follows that $d_0:={\rm dist}(p_0,\widetilde{\partial_H
u}(\tilde \xi))>0$. On account of the upper semicontinuity of
$\widetilde{\partial_H u}$ at $\tilde \xi$, there exists $\delta>0$
such that
$$\widetilde{\partial_H u}(\xi')\subset \widetilde{\partial_H
u}(\tilde \xi)+B_{\mathbb R^{2n}}(0,d_0/4),\ \forall \xi'\in
B_{\mathbb R^{2n}}(\tilde \xi,\delta)\cap \texttt{\rm
Pr}_1(\overline {\Omega_0} \cap H_{\xi_0}).$$ Applying the latter
relation for $\xi'=\tilde\zeta_k$, and taking into account that
$p_k\to 0,$ we obtain that for  $k$ large enough,
$$p_0\in \widetilde{\partial_H
u}(\tilde\zeta_k)+p_k\subset \widetilde{\partial_H u}(\tilde
\xi)+B_{\mathbb R^{2n}}(0,d_0/2),$$ which contradicts the definition
of  $d_0.$ Therefore, $p_0\in \widetilde{\partial_H u}(\tilde \xi).$

We claim that  $\tilde \xi\in \texttt{\rm Pr}_1({\Omega_0}\cap
H_{\xi_0}).$ To see this, we assume by contradiction that  $\tilde
\xi\in \texttt{\rm Pr}_1(\partial{\Omega_0}\cap H_{\xi_0}).$ Then,
$p_0\in \widetilde{\partial_H u}(\tilde \xi)$ is equivalent to $0\in
\mathcal G_1(\tilde \xi)$, which contradicts relation
(\ref{gamma-nulla}). Consequently, $\tilde \xi\in \texttt{\rm
Pr}_1({\Omega_0}\cap H_{\xi_0})$; therefore,
$$p_0\in \widetilde{\partial_H u}(\tilde \xi)={\partial_H u}(\texttt{\rm Pr}_1^{-1}(\tilde \xi))={\partial_H u}(\xi),$$
where $\xi=\texttt{\rm Pr}_1^{-1}(\tilde \xi)\in {\Omega_0}\cap
H_{\xi_0},$ which
concludes the proof. \hfill $\square$\\

\subsection{Comparison principles for $H-$convex functions}

In this subsection we  apply Theorem \ref{n-comparison-uj} to prove
Theorem \ref{comparison-theorem-new} and Theorem \ref{n-comparison
intro}. To do this, we shall compare $H-$convex functions with
specific cone functions, that we will call slicing cones. Some
properties on the horizontal normal mapping of such cones will be
presented in the sequel.

We present in the sequel the construction of this specific cone
function, taking into account that we are in a domain that is
horizontally bounded (but it could be in general, unbounded).

Let $G_0\subset\H^n$ be an open and horizontally bounded set and
$\xi_0\in G_0$ such that $G_0\cap H_{\xi_0}$ is (Euclidean) convex.
Let $c_v<c_b\leq 0$.

For every $\xi \in  H_{\xi_0}$ with $\xi\not=\xi_0,$ we define
$\xi^\partial=\xi^\partial(\xi)$ the unique point in $\partial
G_0\cap H_{\xi_0}$ such that $\xi$ belongs to the horizontal segment
(that is exactly the geodesic in the Carnot-Carath\'eodory metric)
from $\xi_0$ to $\xi^\partial.$ Moreover, for every such $\xi\in
H_{\xi_0}$ with $\xi\not=\xi_0,$ we define $ \lambda^\xi$ as the
unique positive value such that
\begin{equation} \label{lambda-prop}
\xi=\xi_0\circ \delta_{\lambda^\xi}(\xi^{-1}_0\circ \xi^\partial).
\end{equation}
For $\xi= \xi_0$ we set $\lambda^{\xi_0}=0$, we also define
$\xi_{0}^{\partial}$ to be an arbitrary point in $\partial G_{0}\cap
H_{\xi_{0}}$.

Now, for every $\xi\in\H^n,$ we define  $\xi^\perp\in H_{\xi_0}$ to
be the Euclidean orthogonal projection of $\xi$ on the plane
$H_{\xi_0}.$ Finally we define the {\it slicing cone} $V:\mathbb
R^{2n+1}\to\R$  with vertex $(\xi_0,c_v)$ and base $G_0\cap
H_{\xi_0}$ with the value $c_b$ on $\partial G_0\cap H_{\xi_0}$ by
\begin{equation}\label{def v-cone-slide}
V(\xi)=c_v\left(1-\left(1-\frac{c_b}{c_v}\right)\frac{N(\xi_0^{-1}\circ
\xi^\perp)}{N(\xi_0^{-1}\circ (\xi^\perp)^\partial)}\right),\qquad
\xi\in \H^n= \mathbb R^{2n+1}.
\end{equation}
 An easy computation shows that
\begin{equation}\label{v-cone-horizontal}
V(\xi)=c_v\left(1-\left(1-\frac{c_b}{c_v}\right)\lambda^{\xi^\perp}\right),\qquad
\xi\in \H^n.
\end{equation}
Since $\lambda^{\xi^\perp}=\lambda^\xi=1,$ for every $\xi\in
\partial G_0\cap H_{\xi_0},$  we have $V(\xi)=c_b.$

By its definition, the function $V\bigl|_{H_{\xi_0}}$ is Euclidean
convex which implies that $V$ is Euclidean convex and hence
$H-$convex.

\begin{proposition}\label{cone-function-properties}
Let $\Omega\subset \mathbb H^n$ be an open,  horizontally bounded
set, $G_0\subset \Omega$ be an open (Euclidean) convex set,
$\xi_0\in G_0$ and $c_v<c_b\leq 0$. The slicing cone $V:\mathbb
R^{2n+1}\to \mathbb R$ with vertex $(\xi_0,c_v)$ and base $G_0\cap
H_{\xi_0}$ with the value $c_b$ on $\partial G_0\cap H_{\xi_0}$ has
the following properties:
\begin{itemize}
\item[{\rm (i)}]
${B_{\R^{2n}}(0,r_0)}\subset\partial_H V(\xi_0),$ where
$r_0=\frac{c_b-c_v}{{\rm diam}_{H}(G_0\cap H_{\xi_{0}})}; $

\item[{\rm (ii)}]
for every $p\in \texttt{\rm int}(\partial_H V(\xi_0)),$ we have
\begin{equation}\label{big U1-prop-cone-constr}
V(\xi)> V(\xi_0)+ p\cdot ({\rm Pr}_1(\xi)-{\rm Pr}_1(\xi_0)),\qquad
\forall\xi\in \overline {G_0}\cap H_{\xi_0}\setminus \{\xi_0\}.
\end{equation}
\end{itemize}
\end{proposition}

{\it Proof.}
 Let us prove first (i). By definition,   $p\in
\partial_H V(\xi_0)$ is equivalent to the inequality
\begin{equation}\label{dim3-prop}
V(\xi)\ge V(\xi_0)+p\cdot(\texttt{\rm Pr}_1(\xi)-\texttt{\rm
Pr}_1(\xi_0)),\qquad \forall \xi\in G_0\cap H_{\xi_0}.
\end{equation}
We shall use that $V$ on $G_0\cap H_{\xi_0}$ is defined by
\eqref{v-cone-horizontal}, with $\xi^\perp=\xi$. Applying a group
multiplication to the relation  \eqref{lambda-prop} by $\xi_0^{-1}$
from the left and applying the projection map $\texttt{\rm Pr}_1$ to
both sides we obtain
$$\texttt{\rm
Pr}_1(\xi)-\texttt{\rm Pr}_1(\xi_0)=\lambda^\xi(\texttt{\rm
Pr}_1(\xi^\partial)-\texttt{\rm Pr}_1(\xi_0)).$$ Therefore,
\eqref{dim3-prop} is equivalent to the inequality
\begin{equation}\label{dim4-prop}
c_b-c_v \ge p\cdot(\texttt{\rm Pr}_1(\xi^\partial)-\texttt{\rm
Pr}_1(\xi_0)),\qquad \forall \xi\in G_0\cap H_{\xi_0}.
\end{equation}
Since
$$|\texttt{\rm Pr}_1(\xi^\partial)-\texttt{\rm
Pr}_1(\xi_0)|= N(\xi_{0}^{-1}\circ \xi^\partial)\leq {\rm
diam}_{H}(G_0\cap H_{\xi_{0}}),$$by the definition of the number
$r_0>0$  it is easy to see that for all $p\in B_{\mathbb
R^{2n}}(0,r_0),$ relation (\ref{dim4-prop}) holds.

Now, we are going to prove (ii).  Since $\partial_{H}V(\xi_0)$ is
convex and $0\in B_{\mathbb R^{2n}}(0,r_0)\subset\partial_H
V(\xi_0)$ (cf. (i)), $\partial_H V(\xi_0)$ is a star-shaped set with
respect to the origin of $\mathbb R^{2n}.$ Moreover,  $$\texttt{\rm
int}(\partial_H V(\xi_0))=\bigcup\{\alpha p:\alpha\in [0,1),p\in
\partial_H V(\xi_0)\}.$$

Let $\alpha\in (0,1)$ and $p\in\partial_H V(\xi_0)$ be fixed. The
latter relation implies that for every $\beta\in (0,1)$ we have that
$\beta p \in
\partial_{H}V(\xi_0)$, and for every
$\xi\in \overline G_0\cap H_{\xi_0}$,
\begin{equation}\label{proof v2-prop}
V(\xi)\ge V(\xi_0)+ \beta p\cdot(\texttt{\rm Pr}_1(\xi)-\texttt{\rm
Pr}_1(\xi_0)).
 \end{equation}
If $p\cdot(\texttt{\rm Pr}_1(\xi)-\texttt{\rm Pr}_1(\xi_0))>0,$ we
set $\beta=(\alpha+1)/2$ and \eqref{proof v2-prop} implies
\begin{equation}\label{proof v3-prop}
V(\xi)> V(\xi_0)+ \alpha p\cdot(\texttt{\rm Pr}_1(\xi)-\texttt{\rm
Pr}_1(\xi_0)).
 \end{equation}
If $p\cdot(\texttt{\rm Pr}_1(\xi)-\texttt{\rm Pr}_1(\xi_0))<0,$ we
set $\beta=\alpha/2$ and \eqref{proof v2-prop} implies
\begin{equation}\label{proof v4-prop}
V(\xi)> V(\xi_0)+ \alpha p\cdot(\texttt{\rm Pr}_1(\xi)-\texttt{\rm
Pr}_1(\xi_0)).
 \end{equation}
The third possibility is the case when $p\cdot(\texttt{\rm
Pr}_1(\xi)-\texttt{\rm Pr}_1(\xi_0))=0$ for some $\xi \in \overline
{G_0} \cap H_{\xi_0}\setminus \{\xi_0\}$. Since
$V(\xi)=c_v\left(1-\left(1-\frac{c_b}{c_v}\right)\lambda^\xi\right)
> c_v=V(\xi_0)$ we obtain again the inequality
$$
V(\xi)> V(\xi_0)=V(\xi_0)+ \alpha p\cdot(\texttt{\rm
Pr}_1(\xi)-\texttt{\rm Pr}_1(\xi_0)).
$$
 Combining the latter relation with \eqref{proof v3-prop} and
\eqref{proof v4-prop}, we have that  for all $\xi\in \overline
{G_0}\cap H_{\xi_0}\setminus \{\xi_0\},$
$$
V(\xi)> V(\xi_0)+ \alpha p\cdot(\texttt{\rm Pr}_1(\xi)-\texttt{\rm
Pr}_1(\xi_0)),
$$
which concludes the proof.
\hfill $\square$\\

{\it Proof of Theorem \ref{comparison-theorem-new}.} Let $\xi_0\in
\Omega_0$ be fixed.  Without loss of generality,  we may assume that
$u(\xi_0)<v(\xi_0)<0$; otherwise, we subtract a sufficiently large
number from both functions. Let us fix $q\in\partial_H v(\xi_0)$ and
consider the function $U:\overline{\Omega} \to\R$ defined by
$$U(\xi)=u(\xi)-q\cdot ({\rm
Pr}_1(\xi)-{\rm Pr}_1(\xi_0)).$$ Clearly, $U$ is $H-$convex,
$U(\xi_0)=u(\xi_0),$ and
\begin{eqnarray}
U(\xi)&=&u(\xi)-q\cdot ({\rm Pr}_1(\xi)-{\rm
Pr}_1(\xi_0))\nonumber\\
&=&v(\xi)-q\cdot ({\rm Pr}_1(\xi)-{\rm Pr}_1(\xi_0))\nonumber\\
&\ge& v(\xi_0)=u (\xi_0)+m_0,\qquad\forall\xi\in \partial
\Omega_0\cap H_{\xi_0}\label{big U-prop-masik}
\end{eqnarray}
where $m_0=v(\xi_0)-u(\xi_0)>0.$ We notice that for every $\xi\in
\Omega_0$,
\begin{equation}\label{big U00-prop}
\partial_H U(\xi)=\partial_H u(\xi)-q.
\end{equation}

Now let us denote here and in the sequel by $\Omega_0^{\rm conv}$
the Euclidean convex hull of $\Omega_0$ and  we  consider the
slicing cone  $V: \mathbb R^{2n+1}\to\R$
 with vertex $(\xi_0,u(\xi_0))$  and base
$\Omega_0^{\rm conv}\cap H_{\xi_0}$ with the value $v(\xi_0)=u
(\xi_0)+m_0$ on $\partial \Omega_0^{\rm conv}\cap H_{\xi_0};$ see
\eqref{def v-cone-slide}. We know that $V$ is Euclidean convex and
hence $H-$convex.

Since $\Omega_0\subset \Omega_0^{\rm conv}$, from  \eqref{big
U-prop-masik} we have
\begin{equation}\label{big U0-prop}
U(\xi_0)=u(\xi_0)=V(\xi_0)\qquad\texttt{\rm and}\qquad U(\xi)\ge u
(\xi_0)+m_0\ge V(\xi), \quad\forall\xi\in \partial \Omega_0\cap
H_{\xi_0}.
\end{equation}
In addition, by applying Proposition \ref{cone-function-properties}
with $G_0=\Omega_0^{\rm conv}$, $c_b=u(\xi_0)+m_0$ and
$c_v=u(\xi_0)$, and taking into account that  $\partial
\Omega_0\subset \overline{\Omega_0^{\rm conv}},$ we have
\begin{itemize}
\item[(i)]
 ${B_{\R^{2n}}(0,r_{\xi_0})}\subset\partial_H V({\xi_0}),$ where
$r_{\xi_0}=\frac{m_0}{{\rm diam}_{HS}(\Omega_0^{\rm conv})}; $

\item[(ii)] for every $p\in \texttt{\rm int}(\partial_H V({\xi_0})),$ we have
\begin{equation}\label{big U1-prop}
V(\xi)> V({\xi_0})+ p\cdot ({\rm Pr}_1(\xi)-{\rm
Pr}_1({\xi_0})),\qquad \forall\xi\in \partial  \Omega_0\cap
H_{{\xi_0}}.
\end{equation}
\end{itemize}
 Taking into consideration
\eqref{big U0-prop} and (ii) we can apply Theorem
\ref{n-comparison-uj} for the functions $U$ and $V$ on the open
bounded set $\Omega_0\subset \Omega$ to conclude that for any $p\in
\texttt{\rm int}(\partial_H V({\xi_0})),$ we have  $p \in
\partial_H U(\Omega_0\cap H_{{\xi_0}})$. Consequently, one has
 \begin{equation}\label{big U2}
\texttt{\rm int}(\partial_H V({\xi_0}))\subset\partial_H
U(\Omega_0\cap H_{{\xi_0}}).
 \end{equation}
By using (i) and \eqref{big U00-prop} we deduce the following chain
of inclusions:
\begin{equation}\label{non-zero-ball}
0\in{B_{\R^{2n}}\left(0,{r_{\xi_0}}/{2}\right)}\subset\texttt{\rm
int}(\partial_H V({\xi_0}))\subset\partial_H U(\Omega_0\cap
H_{{\xi_0}})=\partial_H u(\Omega_0\cap H_{{\xi_0}})-q.
\end{equation}
In particular, $q\in \partial_H u(\Omega_0\cap H_{{\xi_0}})$, which
concludes the proof. \hfill $\square$\\

The following result is a direct consequence of Theorem
\ref{n-comparison-uj}.

\begin{theorem}\label{cor-strict} Let $\Omega \subset \H^{n}$ be an
open, horizontally bounded and $H-$convex set, $u:\Omega\to\R$ be an
$H-$convex function, and $v:\Omega\to\R$ be a strictly $H-$convex
function. Let ${\Omega_0}\subset \mathbb H^n$ be open such that
$\overline {\Omega_0}\subset\Omega$ and assume that $u\leq v$ in
${\Omega_0}$ and $u=v$ on $\partial {\Omega_0}$. Then
$$\partial_{H}v({\Omega_0}) \subset
\partial_{H}u({\Omega_0}).$$
\end{theorem}

\begin{remark}\rm  The two consequences of Theorem \ref{n-comparison-uj}, i.e. the statements of Theorem
\ref{comparison-theorem-new} and Theorem \ref{cor-strict}, can be
merged once we replace $u<v$ by $u\leq v$ in $\Omega_0$ in the
former, and the strict $H-$convexity by the $H-$convexity in the
latter result. We think that such a general statement is still valid
in our context but the method of Theorem \ref{n-comparison-uj} does
not seem to work. However, Theorem \ref{n-comparison-uj}  is
sufficient to prove the Aleksandrov-type estimate.
\end{remark}

Another consequence of Theorem \ref{n-comparison-uj} is the
Heisenberg comparison principle which corresponds to the Euclidean
one, see Guti\'errez \cite[Theorem 1.4.6]{Gu2001}.

\medskip

 {\it Proof of Theorem \ref{n-comparison intro}.} Without loss of generality, we may assume that $u$ and
 $ v$ are strictly negative in $\overline\Omega$ and that
 $\min_{\xi\in\partial{\Omega}}(v(\xi)-u(\xi))=0.$ Otherwise, we may replace $v$ by $\tilde
v=v+A-\min_{\xi\in\partial{\Omega}}(v(\xi)-u(\xi))$ and $u$ by
$\tilde u=u+A,$ where $A$ is a sufficiently small negative number.

 Suppose that there exists $\xi_0\in\Omega$ such that
 $v(\xi_0)<u(\xi_0)<0.$ Let us fix  $\alpha\in(0,1)$ such that
  $v(\xi_0)< \alpha v(\xi_0)<u(\xi_0)$ and consider the set
$$\Omega_0=\{\xi\in\Omega:\ \alpha v(\xi)<u(\xi)\}.$$
Since $u$ and $v$ are continuous functions on $\Omega$, and
$\xi_0\in \Omega_0,$ it follows that $\Omega_0$ is a non-empty open
set.

We first notice that $\overline{\Omega_0}\subset \Omega.$ Indeed, if
we assume by contradiction that there exists  $\zeta\in \partial
\Omega\cap \overline{\Omega_0}$, then $\alpha v(\zeta)\leq
u(\zeta)$. Since $\min_{\xi\in\partial{\Omega}}(v(\xi)-u(\xi))=0,$
we have that $v(\zeta)\geq u(\zeta),$ a contradiction with the facts
that $\alpha\in (0,1)$ and $u,v$ are strictly negative.

We can apply Theorem \ref{comparison-theorem-new} to functions
$\alpha v < u$ in $\Omega_{0}$ obtaining that $\partial_H
u(\Omega_0)\subset
\partial_H(\alpha v)(\Omega_0)=\alpha\partial_H v(\Omega_0).$ We notice that from the proof of
Theorem \ref{comparison-theorem-new}, by replacing $u$ by $\alpha v$
and $v$ by $u$, respectively,   it also follows that
$\L^{2n}(\partial_{H}v(\Omega_{0}))>0$, see  relation
(\ref{non-zero-ball}). Moreover, by Proposition \ref{usco} one also
has that $\L^{2n}(\partial_{H}v(\Omega_{0}))<+\infty$. Therefore, we
obtain
$$
\L^{2n}(\partial_H u(\Omega_0))\le \alpha^{2n}\L^{2n}(\partial_H
v(\Omega_0))<\L^{2n}(\partial_H
v(\Omega_0)),$$ which contradicts the assumption. \hfill $\square$\\

{\it Proof of Corollary \ref{corollary n-comparison bis convex}.} It
follows directly from Theorem \ref{n-comparison intro}.\hfill $\square$\\

\section{Aleksandrov-type maximum principles}\label{sect-4}

In this section we prove the main result of the paper, i.e., the
Heisenberg version of Aleksandrov's maximum principle in Theorem
\ref{Alexandrov-n}. The proof of Theorem \ref{Alexandrov-n} is based
on a strategy
 following three arguments:

\begin{itemize}
  \item Using the basic comparison principle we shall prove first an
  Aleksandrov-type estimate with respect to {\it horizontal planes}, i.e.,
\begin{equation}\label{alexandrov bis-horizontal-0}
    |u(\xi_0)|^{2n}\le C_n'  \texttt{\rm
dist}_H(\xi_0,\partial\Omega\cap H_{\xi_0})\texttt{\rm
diam}_{H}(\Omega\cap H_{\xi_{0}})^{2n-1}\L^{2n}(\partial_H
u(\Omega\cap H_{\xi_0})),\ \forall \xi_0\in \Omega,
\end{equation}
where $C_n'>0$ depends only on $n$, see Theorem \ref{theo alex}.
Observe that for bounded cylindrical-type domains (which have 'flat
faces' close, but parallel to horizontal planes at a given point)
one may occur that $\texttt{\rm dist}_H(\xi_0,\partial\Omega\cap
H_{\xi_0})\not\rightarrow 0$ in spite of the fact that $\xi_0\to
\partial \Omega$. In such cases the estimate \eqref{alexandrov bis-horizontal-0} is much weaker than the desired \eqref{alexandrov bis}.
The solution to this problem is to compare the values $u(\xi_{0})$
and $u(\zeta)$ where $\zeta\in \Omega$ are close enough  to
$\xi_{0}$ and a better estimate for $\texttt{\rm
dist}_H(\zeta,\partial\Omega\cap H_{\zeta})$ is available.

  \item We establish a  Harnack-type inequality by proving that there exists a constant $C_{1}>1$ such that  if $B_H(\xi_0,3R)\subset \Omega$ for some $\xi_0\in \Omega$ and $R>0$,
 then
$$
    \frac{1}{C_{1}}{u(\xi)}\geq {u(\zeta)} \geq C_{1}{u(\xi)},\quad \forall \xi,\zeta\in B_H(\xi_0,R),
$$
see Theorem \ref{theo harnack-new} in the Appendix. Now, from
(\ref{alexandrov bis-horizontal-0}) and Harnack estimate we have
that
$$ |u(\xi_0)|^{2n}\le C_n''  \mathcal D(\xi_0)\texttt{\rm
diam}_{HS}(\Omega)^{2n-1}\L_{HS}^{2n}(\partial_H u(\Omega)),\
\forall \xi_0\in \Omega,$$ where $C_n''=(C_{1})^{2n}C_n'$ and
$$ \mathcal D(\xi_0)=\min \{ \texttt{\rm dist}_H(\zeta,\partial\Omega\cap H_\zeta) :
\zeta\in \overline {B_H(\xi_0, \texttt{\rm
dist}_H(\xi_0,\partial\Omega)/3)} \}.$$
  \item Finally, by exploiting a typically Heisenberg phenomenon, i.e., the twirling effect of the horizontal planes from one point to another,
 we prove that there is a constant $C_{2}>0$ such that
$$\mathcal D(\xi_0) \leq
C_{2} \texttt{\rm dist}_H(\xi_0,\partial\Omega),\ \ \forall \xi_0\in
\Omega.$$
\end{itemize}

\subsection{Maximum principle on horizontal planes} \label{sect 4.1}

The first step in our strategy consists of the following statement:

\begin{theorem}[Aleksandrov-type maximum principle on horizontal planes]\label{theo alex}

Let $\Omega \subset \H^{n}$ be an open, horizontally bounded and
convex set. If $u:\overline{\Omega} \to \R$ is a continuous
$H-$convex function which verifies $u=0$ on $\partial \Omega$, then
\begin{equation} \label{alexandrov bis-horizontal}
 |u(\xi_0)|^{2n}\le C_n' \texttt{\rm
dist}_H(\xi_0,\partial\Omega\cap H_{\xi_0})\texttt{\rm
diam}_{H}(\Omega\cap H_{\xi_{0}})^{2n-1} \L^{2n}(\partial_H
u(\Omega\cap H_{\xi_0})),\qquad \forall \xi_0\in \Omega,
\end{equation}
where $C_n'>0$ depends only on the dimension $n$.
\end{theorem}

{\it Proof.}  By Proposition \ref{prop-non-positive-H-convex}, we
know that either $u\equiv 0$ on $\overline{\Omega}$, or $u<0$ in
$\Omega$. In the first case, relation (\ref{alexandrov
bis-horizontal}) is trivial; thus, we assume that $u<0$ in $\Omega$.
Let $\xi_0\in\Omega$ be fixed; thus,
 $u(\xi_0)<0.$ The main ingredient of the proof is the
application of Theorem \ref{n-comparison-uj} for an appropriately
constructed comparison function to our function $u$. The proof is
divided into three steps.

{\bf Step 1.} Let $\varepsilon >  0$ be small enough and let
$\Omega_{\varepsilon}$ be an open and convex set (in the Euclidean
sense) such that $\overline{\Omega_{\varepsilon}}\subset \Omega$ and
$\lim_{\varepsilon\to 0^+}\Omega_{\varepsilon}= \Omega$. The
strategy is to prove (in step 2) the Aleksandrov-type
 estimate for the function $u$ restricted to $\Omega_{\varepsilon}$ by means of a  comparison function;
 in step 3, we let $\varepsilon\to 0$. To do this, let us define first the quantity
\begin{equation}\label{megis-kelll}
    \tau_{\xi_0}(\varepsilon) = \min \{ u(\xi): \xi \in \partial \Omega_{\varepsilon} \cap H_{\xi_0}\}.
\end{equation}
Since $u=0$ on $\partial \Omega$ and $u$ is continuous on $\overline
\Omega$, we may consider  $\varepsilon$ so small such that $\xi_{0}
\in \Omega_{\varepsilon}$, and $|\tau_{\xi_0}(\varepsilon)| <
|u(\xi_{0})|/2$.
 Let
  $$t_{\xi_0}(\varepsilon)= 1- \frac{\tau_{\xi_0}(\varepsilon)}{u(\xi_{0})}.$$ Note that $1/2 < t_{\xi_0}(\eps) \leq 1$ and  $t_{\xi_0}(\varepsilon) \to 1$ as $\varepsilon \to 0$.
We shall choose  $v_\varepsilon$ to be the slicing cone $
v_\varepsilon: \R^{2n+1}\to\R$
 with vertex $(\xi_0,u(\xi_0))$
and base $\Omega_{\varepsilon}\cap H_{\xi_0}$ with the value
$\tau_{\xi_0}(\varepsilon)$ on $\partial \Omega_{\varepsilon}\cap
H_{\xi_0}$; see \eqref{def v-cone-slide}. We know $v_\eps$ is
Euclidean convex, then it is $H-$convex.

For further use, let us choose $\xi^-_\varepsilon$ on $\partial
\Omega_{\varepsilon }\cap H_{\xi_0}$ with the property that
$$
 N(\xi^{-1}_0\circ \xi^-_\varepsilon)=\min_{\xi'\in
\partial\Omega_{\varepsilon}\cap H_{\xi_0}}N(\xi^{-1}_0\circ \xi').
$$
Note that the point $\xi^-_\varepsilon$ that realizes the previous
minimum, in general, is not unique. Similarly to
\eqref{lambda-prop}, for every
${\xi}\in\overline{\Omega_{\varepsilon}}\cap H_{\xi_0}$ with
${\xi}\not=\xi_0,$ we define
${\xi}_\eps^\partial={\xi}_\eps^\partial({\xi})$ the unique point in
$\partial\Omega_{\varepsilon}\cap H_{\xi_0}$ such that ${\xi}$
belongs to the horizontal segment  from $\xi_0$ to
${\xi}_\eps^\partial$; let  $\lambda_\eps:= \lambda_\eps^{\xi}$ be
the unique number in $(0,1]$ such that
\begin{equation} \label{lambda}
{\xi}=\xi_0\circ \delta_{\lambda_\eps}(\xi^{-1}_0\circ
{\xi}_\eps^\partial).
\end{equation}
For ${\xi}= \xi_0$ we set $\lambda_\eps^{\xi_0}=0,$ furthermore we
set $\xi_{\epsilon}^{\partial}$ to be an arbitrary point in
$\partial \Omega_{\epsilon}\cap H_{\xi_{0}}$.
 Similarly to (\ref{v-cone-horizontal}), the
restriction of $v_\varepsilon$ to $
\overline{\Omega_{\varepsilon}}\cap H_{\xi_0}$ is explicitly given
by the formula
\begin{equation}\label{dim1bis}
v_\varepsilon({\xi})=u(\xi_0)\left(1-t_{\xi_0}(\varepsilon)\lambda_\eps^{\xi}\right),
\
 {\xi}\in \overline{\Omega_{\varepsilon}}\cap H_{\xi_0}
.\quad
\end{equation}

{\bf Step 2.} On account of (\ref{dim1bis}) and (\ref{megis-kelll})
we observe that
\begin{equation}\label{kell-a-tetelhez}
u(\xi_0)=v_\varepsilon(\xi_0)\qquad {\rm  and}\qquad
u({\xi})\geq\tau_{\xi_0}(\varepsilon)= v_\varepsilon({\xi}),\quad
\forall {\xi}\in
\partial \Omega_\varepsilon\cap H_{\xi_0}.
\end{equation}

 We claim the following properties hold:
\begin{itemize}

\item[(i)]  $B_{\mathbb R^{2n}}(0,r_\varepsilon)\subset\partial_H v_\varepsilon(\xi_0)$ for
 $r_\varepsilon=-\displaystyle t_{\xi_0}(\varepsilon) \frac{u(\xi_0)}{{\rm diam}_{H}(\Omega_\eps\cap H_{\xi_{0}})};$

 \item[(ii)] for every $p\in \texttt{\rm int}(\partial_H v_\eps({\xi_0})),$ we have
\begin{equation}\label{big U1-alex}
v_\eps(\xi)> v_\eps({\xi_0})+ p\cdot ({\rm Pr}_1(\xi)-{\rm
Pr}_1({\xi_0})),\qquad \forall\xi\in \partial  \Omega_\eps\cap
H_{{\xi_0}}.
\end{equation}

\item[(iii)] $p^-_\varepsilon=-u(\xi_0)\displaystyle t_{\xi_0}(\varepsilon) \frac{\texttt{\rm Pr}_1(\xi^-_\varepsilon)-\texttt{\rm Pr}_1(\xi_0)}{|\texttt{\rm Pr}_1(\xi^-_\varepsilon)
-\texttt{\rm Pr}_1(\xi_0)|^2}\in\partial_H v_\varepsilon(\xi_0).$
\end{itemize}
Properties (i) and (ii) follow directly from Proposition
\ref{cone-function-properties}. It remains to prove (iii). To do
that,  $p\in
\partial_H v_\varepsilon(\xi_0)$ is equivalent to the inequality
\begin{equation}\label{dim3}
v_\varepsilon({\xi})\ge v_\varepsilon(\xi_0)+p\cdot(\texttt{\rm
Pr}_1({\xi})-\texttt{\rm Pr}_1(\xi_0)),\qquad \forall {\xi}\in
\Omega_{\varepsilon}\cap H_{\xi_0}.
\end{equation}
By (\ref{lambda}) and (\ref{dim1bis}), the latter inequality reduces
to
\begin{equation}\label{dim4}
-u(\xi_0)t_{\xi_0}(\varepsilon)\ge p\cdot(\texttt{\rm
Pr}_1({\xi}_\eps^\partial)-\texttt{\rm Pr}_1(\xi_0)),\qquad \forall
{\xi}\in \Omega_{\varepsilon}\cap H_{\xi_0}.
\end{equation}
By  inserting $p=p^-_\varepsilon$ in  \eqref{dim4}, we obtain that
$$
(\texttt{\rm Pr}_1(\xi^-_\varepsilon)-\texttt{\rm Pr}_1(\xi_0))
\cdot(\texttt{\rm Pr}_1({\xi}_\eps^\partial)-\texttt{\rm
Pr}_1(\xi_0))\le|\texttt{\rm Pr}_1(\xi^-_\varepsilon) -\texttt{\rm
Pr}_1(\xi_0)|^2,\qquad \forall {\xi}\in\Omega_{\varepsilon}\cap
H_{\xi_0}.$$ From general properties of convex domains, see
Rockafellar \cite{Ro1969}, it follows that the above inequality
holds;  (iii) is proven.

By relation (\ref{kell-a-tetelhez}) and (ii), due to Theorem
\ref{n-comparison-uj}, we have that
\begin{equation}\label{inkluzio-eps}
    \texttt{\rm int}\left(\partial_H
    v_\varepsilon(\xi_0)\right)\subseteq
\partial_H u(\Omega_{\varepsilon}\cap H_{\xi_0}).
\end{equation}

{\bf Step 3.} By (i) and (iii) and since $\partial_H
v_\varepsilon(\xi_0)$ is convex, we have that
\begin{equation}\label{meg-ez-is-kell}
    \{\{p^-_\varepsilon \}\cup B_{\mathbb
R^{2n}}(0,r_\varepsilon)\}^{\rm conv}\subseteq \partial_H
v_\varepsilon(\xi_0).
\end{equation}
Consequently, combining (\ref{meg-ez-is-kell}) and relation
(\ref{inkluzio-eps}), it yields that
$${\rm int}\  \{\{p^-_\varepsilon \}\cup B_{\mathbb
R^{2n}}(0,r_\varepsilon)\}^{\rm conv}\subseteq \partial_H u(
\Omega_\eps\cap H_{\xi_0})\subseteq
\partial_H u(\Omega\cap H_{\xi_0}).$$ Therefore, we have
$$
\L^{2n}(\partial_H u(\Omega\cap H_{\xi_0}))\ge  \L^{2n}\left(
\{\{p^-_\varepsilon \}\cup B_{\mathbb
R^{2n}}(0,r_\varepsilon)\}^{\rm conv}\right) \ge c_n\cdot
|p^-_\varepsilon|r_\varepsilon^{2n-1}
$$
for some constant $c_n>0$ depending only on $n$, i.e., from the
definition of $r_\varepsilon$ and $p_\varepsilon^-$, one has
$$ |u(\xi_0)|^{2n} \leq
C_n' \frac{1}{t_{\xi_0}( \varepsilon)^{2n}} |\texttt{\rm
Pr}_1(\xi^-_\varepsilon)-\texttt{\rm Pr}_1(\xi_0)| {\rm
diam}_{H}(\Omega_\eps \cap H_{\xi_{0}})^{2n-1} \L^{2n}(\partial_H
u(\Omega\cap H_{\xi_0})),$$ with  $C_n'=1/c_n>0$. Since ${\rm
diam}_{H}(\Omega_\eps\cap H_{\xi_{0}})\leq {\rm diam}_{H}(\Omega\cap
H_{\xi_{0}})$ and $\xi^-_\varepsilon\in
\partial\Omega_{\varepsilon}\cap H_{\xi_0}\subset \Omega$, we have that $|\texttt{\rm
Pr}_1(\xi^-_\varepsilon)-\texttt{\rm Pr}_1(\xi_0)| \leq \texttt{\rm
dist}_H(\xi_0,\partial\Omega\cap H_{\xi_0})$  which gives
$$ |u(\xi_0)|^{2n} \leq C_n' \frac{1}{t_{\xi_0}(\varepsilon)^{2n}} \texttt{\rm
dist}_H(\xi_0,\partial\Omega\cap H_{\xi_0})\texttt{\rm
diam}_{H}(\Omega\cap H_{\xi_{0}})^{2n-1} \L^{2n}(\partial_H
u(\Omega\cap H_{\xi_0})).$$ Since $t_{\xi_0}(\varepsilon)\to 1$ as
$\varepsilon \to 0$,  we obtain the desired estimate. The proof is
complete.
 \hfill $\square$

 \begin{corollary}\label{cor-41} Under the same assumptions as in Theorem \ref{theo
 alex}, we have
\begin{equation} \label{alexandrov bis-horizontal-slice}
 |u(\xi_0)|^{2n}\le C_n' \texttt{\rm
dist}_H(\xi_0,\partial\Omega\cap H_{\xi_0})\texttt{\rm
diam}_{HS}(\Omega)^{2n-1} \L_{HS}^{2n}(\partial_H u(\Omega)),\qquad
\forall \xi_0\in \Omega.
\end{equation}
 \end{corollary}

\subsection{Maximum principle in convex domains} \label{sect 4.2}

As we already pointed out at the beginning of the section, it can
happen, that $\texttt{\rm dist}_H(\xi,\partial\Omega\cap
H_{\xi})\not\rightarrow 0$ in spite of the fact that $\xi\to
\partial \Omega$, thus the estimate in (\ref{alexandrov bis-horizontal}) is not enough accurate.
However, by combining Theorem \ref{theo alex} (see also Corollary
\ref{cor-41}) and a Harnack type estimate (see Theorem \ref{theo
harnack-new} in the Appendix), we
  obtain

\begin{theorem} \label{combination}
Let $\Omega \subset \H^{n}$ be an open, horizontally bounded and
convex set. If $u:\overline{\Omega} \to \R$ is a continuous
$H-$convex function such that $u=0$ on $\partial \Omega$, then
\begin{equation} \label{alexandrov bis-horizontal-00000}
 |u(\xi)|^{2n}\le C_n'' \mathcal D(\xi)\texttt{\rm
diam}_{HS}(\Omega)^{2n-1} \L_{HS}^{2n}(\partial_H u(\Omega)),\qquad
\forall \xi\in \Omega,
\end{equation}
where $C_n''>0$ depends only on the dimension $n$, and
$$ \mathcal D(\xi)=\min \{ \texttt{\rm dist}_H(\zeta,\partial\Omega\cap H_\zeta) :
\zeta\in \overline {B_H(\xi, \texttt{\rm
dist}_H(\xi,\partial\Omega)/3)} \}.$$
\end{theorem}

To deduce Theorem \ref{Alexandrov-n} from Theorem \ref{combination}
we need the following geometric result, which exploits the twirling
character of the horizontal planes in the Heisenberg framework.

\begin{proposition} \label{geometric estimate-new}
Let $\Omega\subset\H^n$ be an open,  horizontally bounded and convex
set. Then,
\begin{equation}\label{dist-lemma-dist}
    \mathcal D(\xi)
\leq \left(\frac{\sqrt[4]{97}}{2}+\frac{1}{3}\right) \texttt{\rm
dist}_H(\xi,\partial\Omega),\ \forall \xi\in \Omega.
\end{equation}
\end{proposition}

{\it Proof.} After a left-translation argument, it is enough to
prove inequality (\ref{dist-lemma-dist}) for $\xi=0$. Let $d=
\texttt{\rm dist}_H(0,\partial\Omega)>0$ and fix an element
$\xi_0=(x_0,y_0,t_0)\in
\partial \Omega$ such that $d=d_H(0,\xi_0)$. Since
$\Omega$ is convex, we can fix a supporting hyperplane $\pi_{\xi_0}$
at $\xi_0\in \partial \Omega$ which is represented by
$$\pi_{\xi_0}=\{(x,y,t)\in \mathbb H^n: A\cdot(x-x_0)+B\cdot(y-y_0)+c(t-t_0)=0\},$$
for some $A,B\in \mathbb R^n$ and $c\in \mathbb R.$ For the sake of
notations,  we set $a_k=
(a,...,a)\in \mathbb R^k$ for every $a\in \mathbb R$ and  $k\in
\{1,...,n\}.$

{\bf Case 1.} $A=B=0$.  In this particular case, the horizontal
plane $H_{(0_n,0_n,0)}$ and $\pi_{\xi_0}$ are parallel. Let
$\zeta_0=\left(\left(\frac{d_0}{\sqrt{n}}\right)_n,0_n,0\right)\in
\partial B_H(0,d_0)$ where $d_0=d/3$. Let us denote by $L_0$ the
$(2n-1)-$dimensional plane, which is the intersection of the
horizontal plane $H_{\zeta_0}=\{(x,y,t)\in \mathbb
H^n:t=-2\frac{d_0}{\sqrt{n}}(y_1+...+y_n)\}$ and $\pi_{\xi_0}$. Note
that $\texttt{\rm Pr}_1(L_0)$ is a hyperplane in $\mathbb R^{2n}$
whose equation is given by
\begin{equation}\label{hyper-0}
    y_1+...+y_n+\frac{t_0{\sqrt{n}}}{2{d_0}}=0.
\end{equation}
Since $L_0\subset H_{\zeta_0}$, on the account of equation
(\ref{hyper-0}), we have that
\begin{eqnarray*}
 \texttt{\rm dist}_H(\zeta_0,L_0)  &=& \inf_{\zeta\in L_0} d_H\left(\zeta_0,\zeta\right)=\inf_{\zeta\in
L_0}|\texttt{\rm Pr}_1(\zeta)-\texttt{\rm Pr}_1(\zeta_0)|\\
   &=&\inf_{\tilde \zeta\in
\texttt{\rm Pr}_1(L_0)}|\tilde \zeta-\texttt{\rm
Pr}_1(\zeta_0)|=\frac{\frac{|t_0|\sqrt{n}}{2d_0}}{\sqrt{n}}\\&=&\frac{|t_0|}{2d_0}.
\end{eqnarray*}
First, since $\pi_{\xi_0}$ is a supporting hyperplane at $\xi_0\in
\partial \Omega$ to the convex set  $\Omega$, we have that $$\texttt{\rm
dist}_H(\zeta_0,\partial\Omega\cap H_{\zeta_0})\leq \texttt{\rm
dist}_H(\zeta_0,L_0)=\frac{|t_0|}{2d_0}.$$ On the other hand, since
$d=d_H(0,\xi_0)=N(\xi_0)=N(x_0,y_0,t_0)$, then $|t_0|\leq
d^2=9d_0^2.$ Thus,
$$\mathcal D(0)=\min \left\{
\texttt{\rm dist}_H(\zeta,\partial\Omega\cap H_\zeta) :  \zeta\in
\overline {B_H(0, d_0)} \right\}\leq \texttt{\rm
dist}_H(\zeta_0,\partial\Omega\cap H_{\zeta_0})\leq
\frac{|t_0|}{2d_0}\leq\frac{3}{2}d.$$

{\bf Case 2.} $|A|^2+|B|^2\neq 0$. Clearly, after a normalization,
we may assume that $|A|^2+|B|^2=1.$ Let $\zeta_0=(d_0A,d_0B,0)\in
\partial B_H(0,d_0)$ where $d_0=d/3$ as above.  A simple computation shows
that  the
 plane $\pi_{\xi_0}$ is not parallel to the horizontal plane in $\zeta_0,$
$$H_{\zeta_0}=\left\{(x,y,t):t=2d_0({B}\cdot x-{A}\cdot y)\right\}.$$
Let $L_{AB}=\pi_{\xi_0}\cap H_{\zeta_0}$, which is a
$(2n-1)-$dimensional plane. One has that   $\texttt{\rm
Pr}_1(L_{AB})$ is a hyperplane in $\mathbb R^{2n}$ whose  equation
is obtained after the elimination of $t$ from $\pi_{\xi_0}$ and
$H_{\zeta_0}$, i.e.,
\begin{equation}\label{hyper-AB}
   (A+2cd_0B)\cdot x+(B-2cd_0A)\cdot y-A\cdot x_0-B\cdot y_0-ct_0=0.
\end{equation}
Note that
$$|A+2cd_0B|^2+|B-2cd_0A|^2=|A|^2+|B|^2+4c^2d_0^2=1+4c^2d_0^2>0.$$
Taking into account that $L_{AB}\subset H_{\zeta_0}$, we have that
\begin{eqnarray*}
   \texttt{\rm dist}_H(\zeta_0,L_{AB})&=& \inf_{\zeta\in L_{AB}} d_H(\zeta_0,\zeta)=\inf_{\zeta\in L_{AB}}
   |\texttt{\rm Pr}_1(\zeta)-\texttt{\rm Pr}_1(\zeta_0)|\\
   &=& \inf_{\tilde \zeta\in
\texttt{\rm Pr}_1(L_{AB})}|\tilde \zeta-\texttt{\rm Pr}_1(\zeta_0)|\\
&=&\frac{|d_0(A+2cd_0B)\cdot A+d_0(B-2cd_0A)\cdot B-A\cdot
x_0-B\cdot y_0-ct_0|}{\sqrt{1+4c^2d_0^2}}\\
&=&\frac{|d_0-A\cdot x_0-B\cdot y_0-ct_0|}{\sqrt{1+4c^2d_0^2}}\\
&\leq& d_0+\frac{|A\cdot x_0+B\cdot y_0+ct_0|}{\sqrt{1+4c^2d_0^2}}.
\end{eqnarray*}
By Schwartz inequality and from the fact that $|t_0|\leq
d^2=9d_0^2$, it is clear that
\begin{eqnarray*}
  \frac{|A\cdot x_0+B\cdot y_0+ct_0|}{{\sqrt{1+4c^2d_0^2}}} &\leq&
\sqrt{|x_0|^2+|y_0|^2+\frac{t_0^2}{4d_0^2}}
   \leq  \sqrt[4]{1+\frac{t_0^2}{16d_0^4}} \sqrt[4]{(|x_0|^2+|y_0|^2)^2+t_0^2}  \\
   &\leq& \frac{\sqrt[4]{97}}{2}N(x_0,y_0,t_0)= \frac{\sqrt[4]{97}}{2} d.
\end{eqnarray*}
The rest of the proof is similar to the Case 1. The proof is
concluded.\hfill $\square$\\

{\it Proof of  Theorem \ref{Alexandrov-n}.} It follows from Theorem
\ref{combination} and Proposition \ref{geometric estimate-new}.
\hfill $\square$\\

\medskip
We conclude this section showing that the estimate \eqref{alexandrov
bis} in Theorem \ref{Alexandrov-n} has the natural scaling
invariance property with respect to Heisenberg dilations
$\delta_{\lambda}:$
\begin{remark}\label{scaling invariance}\rm
 Let $\Omega \subset \H^{n}$ and  $u:\overline{\Omega} \to \R$ be as in Theorem \ref{Alexandrov-n}.
Let $\lambda>0$ and $\delta_\lambda\Omega$ be the Heisenberg
dilation of the set $\Omega.$ We define the function
$u^\lambda:\overline{\delta_\lambda\Omega} \to \R$
 by $u^\lambda(\xi)=u(\delta_{\frac{1}{\lambda}}(\xi)).$
Then Theorem \ref{Alexandrov-n} gives that
\begin{equation}\label{new remark}
 |u^\lambda(\xi_1)|^{2n}\leq C_n{\rm dist}_H(\xi_1,\partial (\delta_\lambda \Omega))\
{\rm diam}_{HS}(\delta_\lambda
\Omega)^{2n-1}\L_{HS}^{2n}(\partial_Hu(\delta_\lambda\Omega)),\qquad
\forall \xi_1\in \delta_\lambda\Omega.
\end{equation}
 If we consider $\xi_0=\delta_\lambda(\xi_1)$ and taking into
account that
\begin{itemize}
\item ${\rm diam}_{HS}(\delta_\lambda \Omega)=\lambda {\rm
diam}_{HS}(\Omega),$

\item $p\in\partial_H u^\lambda (\xi)$ if and only if $\lambda p\in
\partial_{H} u(\delta_{\frac{1}{\lambda}}(\xi)),$

\item ${\rm dist}_H(\xi_1,\partial (\delta_\lambda \Omega))=\lambda {\rm dist}_H(\xi_0,\partial \Omega),$

\end{itemize}
we obtain that \eqref{new remark} coincides with \eqref{alexandrov
bis}.
 \end{remark}

\section{Examples: sharpness of the results}\label{sect-5}

In this final section we provide explicit examples showing the
sharpness of our results.

\subsection{Failure of comparison
principles in the absence of convexity}\label{sect-51}

In this subsection we provide an example which shows the failure of
the comparison principle for the horizontal normal mapping in the
absence of the convexity of functions. Let
$$\Omega=\left\{(x,y,t)\in \mathbb H^1:x^2+y^2<1,|t|<1\right\},$$
and $u,v\in \Gamma^\infty( \mathbb H^1)$ be defined by
$$u(x,y,t)=t-(1-t^2)g(x,y), \qquad v(x,y,t)=t,$$ where
$$\left\{
  \begin{array}{lll}
   g\in C^\infty(\mathbb R^2)\ {\rm is\ radial},\ 0\leq g\leq
\frac{1}{4},\\ g>0 \ {\rm on}\
A\left(\frac{1}{4},\frac{3}{4}\right)=: S,\ {\rm and}\ g=0 \ {\rm
on} \ \mathbb R^2\setminus S.
  \end{array}
\right.$$ Here, $A(r,R)\subset \mathbb R^2$ is  the standard open
annulus with center $0$ between the radii $r$ and $R$.

It is clear that $u$ is neither convex nor $H-$convex, while $u=v$
on $\partial \Omega$ and $u\leq v$ in $\Omega$. We shall prove that
\begin{equation}\label{equ-elso-gombok}
    B_{\mathbb R^2}\left(0,{1}/{4}\right)\subset
\partial_H v(\Omega)\setminus\partial_H u(\Omega).
\end{equation}
First of all, since $v$ is regular and $H-$convex, for every
$\xi=(x,y,t)\in \Omega$ one has $$\partial_H
v(\xi)=\{(X_1v(\xi),Y_1v(\xi))\}=\{(2y,-2x)\}.$$ Therefore,
$\partial_H v(\Omega)=B_{\mathbb R^2}(0,2).$

Now, we show that $\partial_Hu(\xi)=\emptyset$ for every
$\xi=(x,y,t)\in \Omega$ with $(x,y)\in B_{\mathbb R^2}(0,1/4)$. By
contradiction, if $p_0\in
\partial_Hu(\xi_0)$ for some $\xi_0=(x_0,y_0,t_0)\in \Omega$ with $(x_0,y_0)\in B_{\mathbb R^2}(0,1/4)$, one has in particular
that
\begin{equation}\label{intersect}
    u(\xi)\geq u(\xi_0)+p_0\cdot ({\rm Pr}_1(\xi)-{\rm Pr}_1(\xi_0)),\qquad
\xi\in \Omega\cap H_{\xi_0}\cap H_{(0,0,t_0)}:=L_0.
\end{equation}
Note that $u(\xi)=u(\xi_0)=t_0$ for every $\xi=(x,y,t)\in L_0$ with
$(x,y)\notin S$; thus, by (\ref{intersect}) it follows that
$p_0\cdot ({\rm Pr}_1(\xi)-{\rm Pr}_1(\xi_0))=0$ for every $\xi\in
L_0$. Now, if we consider $\xi=(x,y,t)\in L_0$ such that $(x,y)\in
S$, then (\ref{intersect}) yields the contradiction
$t_0>t_0-(1-t_0^2)g(x,y)=u(\xi)\geq u(\xi_0)=t_0$. This proves that
$\partial_Hu(\xi_0)=\emptyset.$

Finally, we study $\partial_Hu(\xi)$ for $\xi=(x,y,t)\in \Omega$
such that $(x,y)\notin B_{\mathbb R^2}(0,1/4).$
 Since $u$ is smooth in $\Omega$, if $\partial_Hu(\xi)\not=\emptyset,$ then
 $\partial_Hu(\xi)=\{\nabla_H u(\xi)\}$: hence
 one has
$$X_1u=-(1-t^2)g_x(x,y)+2y(1+2tg(x,y)),\quad
Y_1u=-(1-t^2)g_y(x,y)-2x(1+2tg(x,y)).$$ Since $g$ is radial, we have
$g(x,y)=g(r)$ with $r=\sqrt{x^2+y^2}$, thus for $\xi=(x,y,t)\in
\Omega$, we have
\begin{equation}\label{norm-grad-ex}
   (X_1u(\xi))^2+(Y_1u(\xi))^2=(1-t^2)^2g'(r)^2+4r^2(1+2tg(r))^2.
\end{equation}
Now, for every $\xi=(x,y,t)\in \Omega$ such that $\partial_H
u(\xi)\neq \emptyset$ and $(x,y)\notin B_{\mathbb R^2}(0,1/4)$,
since $0\leq g\leq 1/4$, we have
$$(X_1u(\xi))^2+(Y_1u(\xi))^2\geq \frac{1}{16}.$$
Consequently, $\partial_H u(\Omega)\cap B_{\mathbb
R^2}(0,1/4)=\emptyset,$ which proves the claim.

\begin{remark}\rm We cannot expect even
to have $\mathcal L_{HS}^2(\partial_H v(\Omega))\leq\mathcal
L_{HS}^2(\partial_H u(\Omega))$ for functions $u$ and $v$ with $u=v$
on $\partial \Omega$ and $u\leq v$ in $\Omega$ in the absence of
convexity. Indeed, with respect to the previous example we assume in
addition that $|g'|\leq c$ and $0\leq g\leq c$ for some $c>0$. While
$\mathcal L_{HS}^2(\partial_H v(\Omega))=\mathcal L^2(\partial_H
v(\Omega))=4\pi$, by relations (\ref{equ-elso-gombok}) and
(\ref{norm-grad-ex}) we have $$\mathcal L_{HS}^2(\partial_H
u(\Omega))\leq\mathcal L^2(\partial_H u(\Omega))\leq
\left(c^2+4(1+2c)^2-\frac{1}{16}\right)\pi$$ which is smaller than
$4\pi$ for $c>0$ sufficiently small.
\end{remark}

\subsection{Sharpness of the Aleksandrov-type maximum principle}\label{sect-52}

In this subsection we shall study the sharpness of the
Aleksandrov-type maximum principle for the first Heisenberg group
$\mathbb H^1$. More precisely, under the assumptions of Theorem
\ref{Alexandrov-n}, let us assume that for some $s\geq 1$ we have
$$|u(\xi)|^{2}\leq C_1{\rm dist}_H(\xi,\partial \Omega)^s\ {\rm
diam}_{HS}(\Omega)\L_{HS}^{2}(\partial_Hu(\Omega)),\qquad \forall
\xi\in \Omega.\eqno{(A_s)}$$

\begin{theorem}\label{a-sharp}
 $(A_1)$ is sharp, i.e., the exponent $s$ in $(A_s)$ cannot be
greater than $1.$
\end{theorem}

{\it Proof.} By Theorem \ref{Alexandrov-n}, $(A_1)$ holds for every
horizontally  bounded, open and convex set $\Omega \subset \H^1$ and
every continuous $H-$convex
 function $u:\overline{\Omega} \to \R$
which verifies $u=0$ on $\partial \Omega$.

Let $\eps\in (0,1)$ be arbitrarily fixed. Our claim is proved once
we construct a bounded, convex domain  $\Omega$ and a function
$u:\overline{\Omega} \to \R$ with the above properties such that
$\L_{HS}^{2}(\partial_Hu(\Omega))<\infty$, and
\begin{equation}\label{limsup-hater}
    \sup_{\xi\in \Omega} \frac{|u(\xi)|^{2}}{{\rm dist}_H(\xi,\partial
    \Omega)^{1+\eps}}=+\infty.
\end{equation}
To do this, let us choose $\alpha<1$ and $\beta>1$ such that
\begin{equation}\label{alpha-beta}
\alpha = \frac{2\beta}{4\beta -1} + \frac{\epsilon}{4} < \frac{1}{2}
+ \frac{\epsilon}{2}.
\end{equation}
With these choices of $\alpha$ and $\beta$, we consider the domain
 $${\Omega}_{+}:= \left\{ (x,y,t) \in \mathbb H^1:\ x\in (0,1],\ (y^{2}+t^{2})^{\beta}- x^{\alpha} < 0\right\},$$
and its reflection over the plane $x=1$ defined as
\begin{equation}\label{reflection}
    {\Omega}_{-}:= \left\{ (x,y,t) \in \mathbb H^1:\ (2-x,y,t) \in
{\Omega}_{+} \right\}.
\end{equation}
We shall define the functions $u_\pm: \overline{\Omega}_{\pm} \to
\R$ as
\begin{equation} \label{u1}
 u_{+}(x,y,t):= (y^{2}+t^{2})^{\beta}-x^{\alpha} \quad  \text{and}
\quad u_{-}(x,y,t):= u_{+}(2-x,y,t).
 \end{equation}
 Finally, let $\Omega$ be the open and convex set $\Omega=\Omega_+\cup \Omega_-;$ we define $u:\overline \Omega\to \mathbb R$ by
 \begin{equation} \label{u2}
 u(x,y,t) = u_{\pm}(x,y,t) \qquad  \text{for } \ (x,y,t)\in
\overline{\Omega}_{\pm}.
 \end{equation}
By definition, it is immediate that $u\in C(\overline \Omega)$ is a
convex function such that $u=0$ on $\partial \Omega$ and $u<0$ in
 $\Omega.$ Moreover, $u\in \Gamma^\infty(\texttt{\rm int}(\Omega_+))$ and according to
 Theorem \ref{teorema mappa normale}, for every $\xi=(x,y,t)\in\texttt{\rm int}( \Omega_+)$ we
 have that
 \begin{eqnarray}
\partial_{H}u_{+}(\xi) &=&\left\{\nabla_Hu_+(\xi)\right\}=\left\{(X_1u_+(\xi),Y_1u_+(\xi))\right\}\nonumber\\
&=&\left\{ \left( -\alpha x^{\alpha-1}+ 4\beta y t
(y^{2}+t^{2})^{\beta-1}, 2\beta y  (y^{2}+t^{2})^{\beta-1} - 4\beta
x t (y^{2}+t^{2})^{\beta-1} \right)\right\}.\qquad\label{nabla u+}
\end{eqnarray}
Similarly, for every $\xi=(x,y,t)\in\texttt{\rm int}( \Omega_-)$ we
 have that
\begin{equation}\label{nabla u-}
 \partial_{H}u_{-}(\xi) = \left\{\left( \alpha (2-x)^{\alpha-1}+ 4\beta y t
(y^{2}+t^{2})^{\beta-1}, 2\beta y (y^{2}+t^{2})^{\beta-1} - 4\beta x
t (y^{2}+t^{2})^{\beta-1} \right)\right\}. \end{equation}
 For every
$\xi=(x,y,t)\in\texttt{\rm int}(\Omega_+)$ with $0<x\le
\frac{\alpha}{2\beta}$ we have
\begin{eqnarray*}
 &&|X_1u_+(\xi)|\le \alpha x^{\alpha-1}+2\beta x^\alpha\le 2  \alpha
 x^{\alpha-1},\\
 &&|Y_1u_+(\xi)|\le 2\beta
 x^{\alpha\cdot\frac{\beta-1}{\beta}}|y-2xt|\le  6\beta
 x^{\frac{\alpha}{\beta}\cdot(\beta-\frac{1}{2})}.
\end{eqnarray*}
We deduce that
$$ \partial_{H}u_{+}(\texttt{\rm int}(\Omega_+)) \subseteq A_1\cup A_2,$$
where
$$A_1=X_1u_+ \Bigl(\left[\alpha/(2\beta),1\right],[-1,1],[-1,1]\Bigr)\times
Y_1u_+\Bigl(\left[\alpha/(2\beta),1\right],[-1,1],[-1,1]\Bigr)$$
 and
$$A_2= \left\{ (-v,w): v\in [\gamma, \infty), |w|\leq C v^{(\beta-\frac{1}{2})\cdot \frac{1}{\alpha-1}\cdot \frac{\alpha}{\beta}} \right\},$$
where $\gamma$ and  $C$ are positive constants.
 Clearly, the measure
of $A_1$ is finite while for $A_2$, we have
$$ \L^2(A_2) = C \int_{\gamma}^{\infty}   v^{(\beta-\frac{1}{2})\cdot \frac{1}{\alpha-1}\cdot \frac{\alpha}{\beta}} dv $$
that converges if and only if
 $ \alpha > \frac{2\beta}{4\beta-1} .
$ According to our choice from (\ref{alpha-beta}) the above
condition holds, proving that $\L^{2}(\partial_{H}u_+(\texttt{\rm
int}(\Omega_+))) < \infty$. The fact that
$\L^{2}(\partial_{H}u_-(\texttt{\rm int}(\Omega_-))) < \infty$ works
similarly. Moreover, if $\xi=(x,y,t)\in \Omega_+\cap\Omega_-,$ then
$x=1$ and $\partial_H u(\xi)$ is not a singleton: more precisely,
taking into account \eqref{nabla u+} and \eqref{nabla u-},  we
 have that
$$
\partial_H u(\xi)=\left[-\alpha+4\beta yt(y^2+t^2)^{\beta-1},\alpha+4\beta
yt(y^2+t^2)^{\beta-1}\right]\times Y_1 u_+(\xi)
$$
that implies $\partial_H
u(\xi)\subset[-\alpha-2\beta,\alpha+2\beta]\times[-6\beta,6\beta].$
Therefore, {\small
$$\L_{HS}^{2}(\partial_{H}u(\Omega))\leq\L^{2}(\partial_{H}u(\Omega))=\L^{2}(\partial_{H}u_+(\texttt{\rm
int}(\Omega_+)))+\L^{2}(\partial_{H}u_-(\texttt{\rm
int}(\Omega_-)))+\L^{2}(\partial_{H}u(\Omega_+\cap\Omega_-)) <
 \infty.$$}
Let us choose  $(0,0,0) \in \partial \Omega$ and $\xi=(x,0,0) \in
\Omega$ such that $x \to 0^+$.
 Since
 $\texttt{\rm dist}_{H}(\xi, \partial \Omega) $ is comparable to $x>0$ and $2\alpha<1+\eps$ (cf. (\ref{alpha-beta})), it  follows that
$$\frac{|u(\xi)|^{2}}{{\rm dist}_H(\xi,\partial
    \Omega)^{1+\eps}}=\frac{x^{2\alpha}}{{\rm dist}_H(\xi,\partial
    \Omega)^{1+\eps}}\sim x^{2\alpha-1-\eps}\to +\infty\ \ {\rm as}\ x\to 0^+,$$
concluding the proof of (\ref{limsup-hater}). \hfill $\square$

\begin{remark}\rm  Instead of (\ref{alpha-beta}),  let us choose
 the parameters $\alpha$ and $\beta$ as
$$\alpha = \frac{\beta}{3\beta -1} + \frac{\epsilon}{4}<\frac{1}{3}+\frac{\epsilon}{3},$$
for some $\eps\in (0,1).$ Then, the domain and function introduced
in Theorem \ref{a-sharp} can be used to prove the sharpness of the
Aleksandrov-type maximum principle in the Euclidean case $\mathbb
R^3$ as well (see relation (\ref{elso-Alex}) for $n=3$), i.e.,
$$|u(\xi)|^{3}\leq C_1{\rm dist}(\xi,\partial \Omega)\ {\rm
diam}(\Omega)^2 \L^{2}(\partial u(\Omega)),\qquad \forall \xi\in
\Omega.$$ The details are left as an exercise to the interested
reader.
\end{remark}

\subsection{Horizontal Monge-Amp\`ere operator versus horizontal normal
mapping}\label{sect-53}

Let $\Omega\subset \mathbb H^1$ be an open, bounded and convex set.
We consider the horizontal Monge-Amp\`ere operator
\begin{equation}\label{sm_a}
   \mathcal S_{ma}(u)(\xi)= {\rm det} [{\rm
 Hess}_H(u)(\xi)]^*+12(Tu(\xi))^2,
\end{equation}
where $u\in C^2(\Omega)$ and $[{\rm
 Hess}_H(u)(\xi)]^*$ is the symmetrized horizontal Hessian: $$[{\rm
 Hess}_H(u)(\xi)]^*=\left[
  \begin{array}{cc}
    X_1^2u & (X_1Y_1u+Y_1X_1u)/2 \\
    (X_1Y_1u+Y_1X_1u)/2 & Y_1^2u \\
  \end{array}
\right](\xi),\ \xi\in \Omega.$$

Having in our mind relation (\ref{monge-ampere-measure}) from the
Euclidean case, we are interested to study the connection between
the quantities $\int_\Omega \mathcal S_{ma}(u)(\xi)d\xi$ and
$\L^2(\partial_H u(\Omega))$ (or $\L_{HS}^2(\partial_H u(\Omega))$ )
whenever $u\in C^2(\Omega)$ is an $H-$convex function. Some initial
information are available as follows:
\begin{itemize}
  \item In \cite{CaPi2011} the authors prove that
$$
 \int_{\partial_H u(\Omega)}{\cal S}^2_\H\bigl(\{\xi\in
\Omega:\,\nabla_Hu(\xi)=v\} \bigr) \, d v=\int_\Omega  \left({\rm
det} [{\rm
 Hess}_H(u)(\xi)]^*+4(Tu(\xi))^2\right) d \xi,$$ where  ${\cal S}^2_\H$ denotes the
2--dimensional spherical Hausdorff measure. Note that if $u\in
\Gamma^2(\Omega)$ is $H-$convex, the matrix $[{\rm
 Hess}_H(u)(\xi)]^*$ is positive semi-definite for every $\xi\in \Omega$ (see Danielli, Garofalo and Nhieu \cite{DaGaNh2003}), thus the latter
 integral and  $\int_\Omega \mathcal S_{ma}(u)(\xi)d\xi$ are
 comparable.
  \item By the oscillation estimate of Guti\'errez and
Montanari \cite[Theorem 1.4]{GuMo-CommPDE}, we know that for any
compact domain $A\subset \Omega$ there exists a constant
$C=C(A,\Omega)>0$ such that $$\int_A \mathcal S_{ma}(u)(\xi)d\xi
\leq C (\sup_\Omega u-\inf_\Omega u)^2$$ for every $H-$convex
function $u\in C^2(\Omega)$. By combining this result
 with our Aleksandrov-type maximum principle in \eqref{alexandrov bis}, one
has that for every compact set $A\subset \Omega$ and  for every
 $H-$convex function $u\in C(\overline \Omega)\cap C^2(\Omega)$ with $u=0$ on
$\partial \Omega$,
 $$\int_A \mathcal
S_{ma}(u)(\xi)d\xi \leq C_1C {\rm
diam}_{HS}(\Omega)\L_{HS}^2(\partial_H u(\Omega)).$$
\end{itemize}
Clearly, if $\int_\Omega \mathcal S_{ma}(u)(\xi)d\xi$ were
comparable to $\L_{HS}^2(\partial_H u(\Omega))$, then our
Aleksandrov-type maximum principle would provide an estimate of the
form
$$|u(\xi_0)|^{2} \leq C \texttt{\rm dist}_H(\xi_0,\partial \Omega) {\rm diam}_{HS}(\Omega)\int_{\Omega}{\cal S}_{ma}(u)(\xi)d \xi,\qquad \xi_0\in \Omega.$$
Unfortunately, this turns out to be  only a wishful thinking as
shown by the following:


\begin{proposition}\label{prop-MA} There exists an open, bounded and convex set $\Omega\subset \mathbb
H^1$, and an  $H-$convex function $u:\overline \Omega\to \mathbb R$
with $u\in C(\overline \Omega)\cap C^2(\Omega)$, $u=0$ on $\partial
\Omega$,  such that
\begin{itemize}
  \item[{\rm (i)}] $\L_{HS}^2(\partial_H u(\Omega))=\infty ;$
  \item[{\rm (ii)}] $\int_{\Omega}{\cal S}_{ma}(u)(\xi)d
\xi< \infty.$
\end{itemize}
\end{proposition}

{\it Proof.} The construction is similar to \eqref{u1} and
\eqref{u2}. More precisely, let us consider $\beta
>1$ with
$$ \frac{1}{2}<\alpha\le  \frac{2\beta}{4\beta-1},$$
the new domain
 $${\Omega}_{+}:= \left\{ \xi=(x,y,t) \in \mathbb H^1:\ x\in (0,1],\ (y^{2}+t^{2})^{\beta}- x^{\alpha}+\frac{\alpha}{2}x^2 < 0\right\},$$
and its  reflection ${\Omega}_{-}$ over the plane $x=1$ defined as
in \eqref{reflection}. The functions $u_{\pm}:
\overline{\Omega}_{\pm} \to \R$ are defined as
$$ u_{+}(x,y,t):= (y^{2}+t^{2})^{\beta}-x^{\alpha}
+\frac{\alpha}{2}x^2 \quad  \text{and} \quad u_{-}(x,y,t):=
u_{+}(2-x,y,t).
$$
Let $\Omega=\Omega_+\cup \Omega_-$, which is an open and convex set;
 we define $u:\overline \Omega\to \mathbb R$ in the same way as in
 \eqref{u2}.
It is a straighforward computation to see that $u\in C(\overline
\Omega)$, and $u\in C^2( \Omega)$ since
$$\frac{\partial u_+}{\partial x}(1,y,t)=\frac{\partial u_-}{\partial x}(1,y,t)=0\qquad\texttt{\rm
and}\qquad \frac{\partial^2 u_+}{\partial
x^2}(1,y,t)=\frac{\partial^2 u_-}{\partial
x^2}(1,y,t)=-\alpha^2+2\alpha.$$
 Moreover,  $u$
is a convex function on $\Omega$ such that $u=0$ on $\partial
\Omega$ and $u<0$ in
 $\Omega.$

(i) First of all, note that
$$\L_{HS}^{2}(\partial_Hu(\Omega))\geq \limsup_{k\to \infty}
\L^{2}\left(\partial_{H}u\left(A_+\cap
H_{(\frac{1}{2k},0,0)}\right)\right),$$
 where $ A_+=\left\{\xi=(x,y,t)\in\texttt{\rm int}(\Omega_+):\ y\ge
0)\right\}. $ Since $u_+$ is regular and $H-$convex in $\texttt{\rm
int}(\Omega_+)$, we have
 \begin{eqnarray*}
\partial_{H}u(\xi) &=&\left\{\nabla_Hu_+(\xi)\right\}=\left\{(X_1u_+(\xi),Y_1u_+(\xi))\right\}\\
&=&\left\{ \left( -\alpha x^{\alpha-1}+ \alpha x+4\beta y t
(y^{2}+t^{2})^{\beta-1}, 2\beta (y-2xt)  (y^{2}+t^{2})^{\beta-1}
\right)\right\},\ \xi\in \texttt{\rm int}(\Omega_+).
\end{eqnarray*}
Therefore, for every $k\geq 1$, one has { \begin{eqnarray*}
S_k&:=&\partial_H u\left(A_+\cap
H_{(\frac{1}{2k},0,0)}\right)\\
 &=&\left\{\left( \alpha x(1-
x^{\alpha-2})-\frac{4\beta}{k}y^{2\beta}\left(1+\frac{1}{k^2}\right)^{\beta-1},
2\beta
y^{2\beta-1}\left(1+\frac{2x}{k}\right)\left(1+\frac{1}{k^2}\right)^{\beta-1}
\right):\right.\\&&\ 0<x<1,\ 0\le y<
\left.\left(x^\alpha -\frac{\alpha}{2}x^2\right)^{\frac{1}{2\beta}}\left(1+\frac{1}{k^2}\right)^{-1/2}\right\}\\
&\supset & \left\{\left( \alpha x(1-
x^{\alpha-2})-\frac{4\beta}{k}y^{2\beta}\left(1+\frac{1}{k^2}\right)^{\beta-1},
2\beta
y^{2\beta-1}\left(1+\frac{2x}{k}\right)\left(1+\frac{1}{k^2}\right)^{\beta-1}
\right):\right.\\&&\ 0<x<1,\ 0\le y<
\left.\left(\frac{x^\alpha}{2^{\beta+1}} \right)^{\frac{1}{2\beta}}\right\}\\
&\supset & \left\{\left( \alpha x(1-
x^{\alpha-2})-\frac{4\beta}{k}y^{2\beta}\left(1+\frac{1}{k^2}\right)^{\beta-1},
2\beta
y^{2\beta-1}\left(1+\frac{2x}{k}\right)\left(1+\frac{1}{k^2}\right)^{\beta-1}
\right):\right.\\&&\ \left.0 \leq y \leq
2^{-\frac{\beta+1}{2\beta}},\
(2^{\beta+1}y^{2\beta})^\frac{1}{\alpha}< x < 1 \right\}.
\end{eqnarray*}}
By the Fatou lemma, we have that
$$\liminf_{k\to \infty}
\L^{2}(S_k)\geq \mathcal L^2(S),$$ where
$$S=\left\{\left(\alpha x(1-
x^{\alpha-2}), 2\beta y^{2\beta-1}\right): 0 \leq y \leq
2^{-\frac{\beta+1}{2\beta}},\
(2^{\beta+1}y^{2\beta})^\frac{1}{\alpha}< x < 1\right\}.$$ On the
other hand, we have that
 $$\L^2(S)\ge \int_{0}^{\gamma}
\alpha
2^{\frac{\beta+1}{\alpha}}\left(\frac{s}{2\beta}\right)^{\frac{2\beta}{(2\beta-1)\alpha}}
\left(- 1+
2^{\frac{(\alpha-2)(\beta+1)}{\alpha}}\left(\frac{s}{2\beta}\right)^{\frac{2\beta(\alpha-2)}{(2\beta-1)\alpha}}
\right)
 ds,
$$ where $\gamma$ is a positive constant depending only on $\beta.$ The latter integral is $+\infty$ since
 $ \alpha \leq \frac{2\beta}{4\beta-1} .$

(ii) By symmetry, it is enough to prove the claim for $u_{+}$.
 Since $\displaystyle \int_{\Omega_{+}}(T u_{+})^{2} d \xi <
 \infty,$ by \eqref{sm_a} we
only need to consider the integral $\displaystyle \int_{\Omega_{+}}
\texttt{\rm det} [{\rm Hess}_H(u_+)(\xi)]^*
 d \xi$
which is clearly finite  if $\displaystyle\int_{\Omega_{+}}
(X_1^{2}u_{+}Y_1^{2}u_{+})(\xi) d\xi < \infty$. The singular term in
the integral is coming from
$$X_1^{2}u_{+}(x,y,t) = -\alpha(\alpha-1)x^{\alpha-2} +\alpha+ 8\beta y^{2} \frac{\partial}{\partial t}\left(t(y^{2}+t^{2})^{\beta-1}\right). $$
Calculating the term $Y^{2}_1u_{+},$ since $0<x<1,$ we obtain
$$|Y_1^{2}u_{+}(x,y,t)|\le C(y^2+t^2)^{\beta-1},$$
 for some constant $C=C(\beta)>0$.
Using integration in  polar coordinates in the $(y,t)-$plane,  we
have
$$\int_{\Omega_{+}} |X_1^{2}u_{+}Y_1^{2}u_{+}|(\xi) d\xi  \le C'
\int_{0}^{1} x^{\alpha-2} \int_{0}^{x^{\frac{\alpha}{2\beta}}}
r^{2\beta-1} d r dx+C'  =\frac{C'}{2\beta} \int_{0}^{1}x^{2\alpha-2}
d x+C',
$$
for some constant $C'=C'(\alpha,\beta)>0.$
Since $\displaystyle \alpha > \frac{1}{2}$,  the above integral
converges.
 \hfill $\square$

\begin{remark}\rm
Unlike in the Euclidean case (see relation
(\ref{monge-ampere-measure}) versus Proposition \ref{prop-MA}),  the
horizontal normal mapping does not play the same role as the Euclidean
normal mapping in the study of the Monge-Amp\`ere equation via the
operator $\mathcal S_{ma}$ given by (\ref{sm_a}). Furthermore, if
$\Omega\subset \mathbb H^n$ is an open, bounded and convex set, and
$u:\overline \Omega\to \mathbb R$ is a continuous $H-$convex
function, we may consider for every $E\subset \Omega$ the function
$\nu_u(E)=\L^{2n}_{HS}(\partial_H u (E)),$ which is a natural
candidate for the Monge-Amp\`ere measure in the Heisenberg setting.
This defines  an outer measure, however $\nu_u$ is not a Borel
 measure in general. Indeed, let
$\Omega\subset \mathbb H^1$ be the cylinder introduced in
\S\ref{sect-51} and let $D_i=\{(x,y,t)\in \Omega:t=t_i\}$, $i\in
\{1,2\},$ be two discs with $-1<t_1<t_2<1.$ If $u(x,y,t)=t$, then
$\nu_u(D_1\cup D_2) = \nu_u(D_1) = \nu_u(D_2) = 4\pi$, i.e., the
additivity on Borel sets of $\nu_u$ fails.
\end{remark}

\section{Appendix}

\subsection{Degree theory for set-valued
maps}\label{appendix-1}

We recall some facts from the degree theory for upper semicontinuous
set-valued maps, see Hu and Papageorgiou \cite{Hu-Papa}. Note that
the degree theory developed in \cite{Hu-Papa} is also valid for
infinite-dimensional spaces, which is a generalization of the
Brouwer, Browder and Leray-Schauder degree theories.  In our
context, it is enough to consider the finite-dimensional version.

Let us start with the definition of  Brouwer degree $\texttt{\rm
deg}_B$ for a continuous function:

\begin{theorem} {\rm (see \cite{Lloyd1978})} \label{degree-Browder}
Let $$M=\left\{ (f,U,y):\ U \subset \R^{n}\ \textrm{open and
bounded},\ f \in C(\overline{U},\R^{n}),\ y \in \R^{n} \backslash
f(\partial U) \right\}.$$ There exists a function, called the {\rm
Brouwer degree}, $\texttt{\rm deg}_B:M \rightarrow \mathbb{Z}$, that
satisfies the following properties:
\begin{itemize}
 \item if $\texttt{\rm deg}_B(f, U, y) \neq 0$, then there exists $x \in U$ such that $f(x)=y$;
 \item $\texttt{\rm deg}_B(Id, U, y)=1$ if $y \in U$;
 \item if $\mathcal F: [0,1] \times \overline{U} \rightarrow \R^{n}$ is a homotopy
 such that $y \in \R^{n} \backslash \mathcal F([0,1] \times \partial U)$,
then $t \mapsto \texttt{\rm deg}_B(\mathcal F(t,\cdot), U, y)$ is
constant:
 \item $\texttt{\rm deg}_B(f, U, y) = \texttt{\rm deg}_B(f-y, U, 0)$.
\end{itemize}
\end{theorem}

\noindent In order to work with the degree of set-valued maps, we
need the following notion.


\begin{definition} {\rm (see \cite[Definition 3]{Hu-Papa})} \label{P-class}
Let $X$ be a finite-dimensional normed space  and $U\subset X$ be an
open bounded set. A set-valued map $F:\overline U\to 2^{X}\setminus
\{\emptyset\}$  is said to belong to the {\rm class (P)} if:
\begin{itemize}
  \item[{\rm (i)}] it maps bounded sets into relatively compact sets;
  \item[{\rm (ii)}] for every $x\in \overline U$, $F(x)$ is closed and convex in
  $X$;
  \item[{\rm (iii)}] $F$ is upper semicontinuous on $\overline U$.
\end{itemize}
\end{definition}

\noindent
 A parameter-depending version of Definition \ref{P-class} reads as
 follows, which will be used to exploit homotopy properties of certain
 set-valued maps.

\begin{definition} {\rm (see \cite[Definition 9]{Hu-Papa})} \label{P-class-homo}
Let $X$ be a finite-dimensional normed space  and $U\subset X$ be an
open bounded set. A one-parameter family of set-valued maps
$\mathcal F_\lambda:\overline U\to 2^{X}\setminus \{\emptyset\}$,
$\lambda\in [0,1]$ is said to be a {\rm homotopy of class (P)}  if:
\begin{itemize}
  \item[{\rm (i)}] $\overline{\{{\cup \mathcal F_\lambda(x):(\lambda,x)\in [0,1]\times \overline U\}}}$ is compact in $X$;
  \item[{\rm (ii)}] for every $(\lambda,x)\in [0,1]\times \overline U$,  $\mathcal F_\lambda(x)$ is closed and convex in
  $X$;
  \item[{\rm (iii)}] $(\lambda,x)\mapsto \mathcal F_\lambda(x)$ is upper semicontinuous from $[0,1]\times \overline U$ into $2^{X}\setminus \{\emptyset\}$.
\end{itemize}
\end{definition}

\noindent For the set-valued degree of upper semicontinuous
set-valued map certain selectors are needed:

\begin{proposition} {\rm (see \cite{Cellina})} \label{prop-cellina} If $X,V$ are Banach spaces, $U\subset X$ is an open bounded set and $F:\overline U\to 2^V\setminus
\{\emptyset\}$ is an upper semicontinuous set-valued map with closed
and convex values then for every $\varepsilon>0$ there exists a
\texttt{\rm continuous approximate selector}
$f_\varepsilon:\overline U\to V$ such that
$$f_\varepsilon(y)\in F((y+B_X(0,\varepsilon))\cap \overline U)+B_V(0,\varepsilon),\ \ \forall y\in \overline U.$$
\end{proposition}

\noindent The next result is a set-valued version of Theorem
\ref{degree-Browder} and it plays a fundamental role in our degree
theoretical argument from Section \ref{sect-3}.

\begin{theorem} {\rm (see \cite[Definition 11 and Theorem 12]{Hu-Papa}\label{Hu-Papa-theorem})} Let $X$ be a finite-dimensional normed space.
Let
$$  M_{SV}=\left\{(F,U,y):
  \begin{array}{lll}
   U \subset X\ \textrm{open and bounded}, \\
   F:\overline U\to 2^{X}\setminus \{\emptyset\}\ {belongs\ to\ the\
class\ {\rm (P)}}, \ y \in X \backslash F(\partial U)
  \end{array}
\right\}.
$$
 There
exists a function, called as  {\rm set-valued degree function},
$\texttt{\rm deg}_{SV}:M_{SV} \rightarrow \mathbb{Z}$, that is
defined as the common value
$$\texttt{\rm deg}_{SV}(F,U,y)=\texttt{\rm deg}_B(f_\varepsilon,U,y)$$ for every small
$\varepsilon>0$, where  $f_\varepsilon$ comes from Proposition {\rm
\ref{prop-cellina}}. The function $\texttt{\rm deg}_{SV}$ verifies
the properties of
\begin{itemize}
  \item {\rm normalization}: $\texttt{\rm deg}_{SV}(Id,U,y)=\texttt{\rm deg}_{B}(Id,U,y)=1$ for all $y\in U;$
  \item {\rm additivity on domain}: If $U_1,U_2\subset U$ are
  disjoint open sets and $y\notin F(\overline U\setminus (U_1\cup
  U_2)),$ then
  $$\texttt{\rm deg}_{SV}(F,U,y)=\texttt{\rm deg}_{SV}(F,U_1,y)+\texttt{\rm deg}_{SV}(F,U_2,y);$$
  \vspace{-0.9cm}
  \item {\rm homotopy invariance}:  if $\mathcal F_\lambda: \overline U\to
2^X$ is a homotopy of class {\rm (P)} and $\gamma:[0,1]\to X$ is
such that $\gamma(\lambda)\notin \mathcal F_\lambda(\partial U)$ for
all $\lambda\in [0,1]$, then $\texttt{\rm deg}_{SV}(\mathcal
F_\lambda,U,\gamma(\lambda))$ is independent of $\lambda\in [0,1].$
\end{itemize}
\end{theorem}

\subsection{Quantitative Harnack-type inequality for $H-$convex functions}\label{appendix-2}

\begin{lemma}\label{two-sided-Harnack-lemma}  Let $\Omega$ be an open convex domain such that $B_H(0,cR)\subset \Omega$ for some constants $c,R>0.$  Let  $u:\overline{\Omega}\to\R$ be
an $H-$convex function with $u\leq 0$ in $\Omega$. Let
$\xi_1,\xi_2\in B_H(0,cR)$ with $\xi_2\in H_{\xi_1}$ and some
constants $c_1,c_2\geq 0$ and $c_3>0$ such that
$$N(\xi_1)\leq c_1R;\ N(\xi_2)\leq c_2R,\ d_H(\xi_1,\xi_2)\leq c_3R
$$ and
$$c_1+c_3<c;\ c_2+c_3<c.$$
Then $$\frac{c-c_1-c_3}{c-c_1}{u(\xi_1)}\geq {u(\xi_2)} \geq
\frac{c-c_2}{c-c_2-c_3}{u(\xi_1)}.$$
\end{lemma}

{\it Proof.} The idea of the proof is close to Lemma 5.2 from
Guti\'errez and Montanari \cite{GuMo-CommPDE}.  Let
$\xi'_\lambda=\xi_1\circ \delta_\lambda(\xi_1^{-1}\circ \xi_2)\in
H_{\xi_1}$ for $\lambda>0$. If $\xi'_\lambda\in \partial B_H(0,cR)$,
then we have that
\begin{eqnarray*}
  cR &=& N(\xi'_\lambda)= N(\xi_1\circ \delta_\lambda(\xi_1^{-1}\circ \xi_2))\\
   &\leq& N(\xi_1)+N(\delta_\lambda(\xi_1^{-1}\circ \xi_2)) \\
   &=&N(\xi_1)+\lambda N(\xi_1^{-1}\circ \xi_2)\\
   &\leq&c_1R+\lambda c_3R.
\end{eqnarray*}
Therefore,
$$\lambda\geq \frac{c-c_1}{c_3}> 1.$$
Now, the relation $\xi'_\lambda=\xi_1\circ
\delta_\lambda(\xi_1^{-1}\circ \xi_2)$ can be written into the form
$\xi_2=\xi_1\circ \delta_{1/\lambda}(\xi_1^{-1}\circ \xi_\lambda')$.
The $H-$convexity of $u$ and the fact that $u\leq 0$ yields that
$$u(\xi_2)\leq \left(1-\frac{1}{\lambda}\right)u(\xi_1)+\frac{1}{\lambda} u(\xi_\lambda')\leq \left(1-\frac{1}{\lambda}\right)u(\xi_1).$$
Consequently, $${u(\xi_2)}\leq
\left(1-\frac{1}{\lambda}\right){u(\xi_1)}\leq
\left(1-\frac{c_3}{c-c_1}\right)u(\xi_1)=\frac{c-c_1-c_3}{c-c_1}u(\xi_1).$$
Now, changing the roles of $\xi_1$ and $\xi_2$, by taking into
account that $\xi_2\in H_{\xi_1}$ (thus, $\xi_1\in H_{\xi_2}$), we
obtain in a similar manner that
$${u(\xi_1)}\leq \frac{c-c_2-c_3}{c-c_2}{u(\xi_2)},$$
which ends the proof. \hfill $\square$\\

\begin{theorem}[Harnack-type inequality] \label{theo harnack-new}
Let $\Omega \subset \H^{n}$ be an open, horizontally bounded and
convex set. If $u:\overline{\Omega}\to\R$ is an $H-$convex function
with $u=0$ on $\partial\Omega$, and $B_H(\xi_0,3R)\subset \Omega$
for some $\xi_0\in \Omega$ and $R>0$, then
\begin{equation} \label{harnack-new}
\frac{1}{31}{u(\xi)}  \geq {u(\zeta)}\geq 31{u(\xi)}, \qquad \forall
\xi,\zeta \in B_H(\xi_0,R).
\end{equation}
\end{theorem}

{\it Proof.} The proof is similar to Guti\'errez and Montanari
\cite[Proposition  5.3]{GuMo-CommPDE}. After a left-translation by
$\xi_0^{-1}$, it is enough to prove (\ref{harnack-new}) for every
$\xi,\zeta\in B_H(0,R)$.

By the first part of Proposition \ref{prop-non-positive-H-convex}
one has that $u\leq 0$ on $\overline\Omega$. Let us fix
$\xi=(x_0,y_0,t_0)\in B_H(0,R)$ arbitrarily, i.e., $N(\xi)\leq R,$
with $x_0=(x^0_1,...,x^0_n)\in \mathbb R^n$ and
$y_0=(y^0_1,...,y^0_n)\in \mathbb R^n$. In particular, we have that
$\sqrt{|t_0|}\leq R.$ For simplicity, we assume that $t_0\geq 0$
(the case $t_0<0$ works similarly).

{\bf Step 1.} Let $\xi_1=\exp\left(-\sum_{j=1}^n
(x_j^0X_j+y_j^0Y_j)\right)\circ\xi=(0_n,0_n,t_0)\in H_{\xi}$. It is
clear that
$$N(\xi)\leq R;\ N(\xi_1)=\sqrt{t_0}\leq R;\ d_H(\xi,\xi_1)=\sqrt{|x_0|^2+|y_0|^2}\leq R.$$
Thus, we may apply Lemma \ref{two-sided-Harnack-lemma} with
$c_1=c_2=c_3=1$ and $c=3$, obtaining
$$\frac{1}{2}{u(\xi)}\geq {u(\xi_1)} \geq
2{u(\xi)}.$$

{\bf Step 2.} Let $\xi_2=\exp\left(\sigma
\sum_{j=1}^nX_j\right)\circ\xi_1=(\sigma_n,0_n,t_0)\in H_{\xi_1}$,
where
$$\sigma=\frac{\sqrt{t_0}}{2\sqrt{n}}.$$ Note that
$$N(\xi_1)\leq R;\ N(\xi_2)=(n^2\sigma^4+t_0^2)^\frac{1}{4}=17^\frac{1}{4}\sigma\sqrt{n}\leq\frac{ 17^\frac{1}{4}}{2}R;\ d_H(\xi_1,\xi_2)=\sigma\sqrt{n}=\frac{\sqrt{t_0}}{2}\leq \frac{R}{2}.$$
Therefore, we  apply Lemma \ref{two-sided-Harnack-lemma} with
$c_1=1,$ $c_2=\frac{ 17^\frac{1}{4}}{2},$ $c_3=\frac{1}{2}$ and
$c=3$, obtaining
$$\frac{\frac{3}{2}}{2}{u(\xi_1)}\geq {u(\xi_2)} \geq
\frac{3-\frac{ 17^\frac{1}{4}}{2}}{\frac{5}{2}-\frac{
17^\frac{1}{4}}{2}}{u(\xi_1)}.$$

{\bf Step 3.} Let $\xi_3=\exp\left(\sigma
\sum_{j=1}^nY_j\right)\circ\xi_2=(\sigma_n,\sigma_n,t_0-2\sigma^2n)\in
H_{\xi_2}$. Note that
$$N(\xi_2)\leq \frac{ 17^\frac{1}{4}}{2}R;\ N(\xi_3)=\sqrt[4]{4\sigma^4n^2+(t_0-2\sigma^2n)^2}=8^\frac{1}{4}\sigma\sqrt{n}\leq\frac{ 8^\frac{1}{4}}{2}R;\ d_H(\xi_2,\xi_3)=\sigma\sqrt{n}\leq \frac{R}{2}.$$
Now, we  apply Lemma \ref{two-sided-Harnack-lemma} with $c_1=\frac{
17^\frac{1}{4}}{2},$ $c_2=\frac{ 8^\frac{1}{4}}{2},$
$c_3=\frac{1}{2}$ and $c=3$, obtaining
$$\frac{\frac{5}{2}-\frac{ 17^\frac{1}{4}}{2}}{3-\frac{
17^\frac{1}{4}}{2}}{u(\xi_2)}\geq {u(\xi_3)} \geq \frac{3-\frac{
8^\frac{1}{4}}{2}}{\frac{5}{2}-\frac{
8^\frac{1}{4}}{2}}{u(\xi_2)}.$$

{\bf Step 4.} Let $\xi_4=\exp\left(-\sigma
\sum_{j=1}^nX_j\right)\circ\xi_3=(0_n,\sigma_n,t_0-4\sigma^2n)=(0_n,\sigma_n,0)\in
H_{\xi_3}$. Note that
$$N(\xi_3)\leq \frac{ 8^\frac{1}{4}}{2}R;\ N(\xi_4)=\sigma\sqrt{n}\leq\frac{R}{2};\ d_H(\xi_3,\xi_4)=\sigma\sqrt{n}\leq \frac{R}{2}.$$
We  apply Lemma \ref{two-sided-Harnack-lemma} with $c_1=\frac{
8^\frac{1}{4}}{2},$ $c_2=c_3=\frac{1}{2}$ and $c=3$, obtaining
$$\frac{\frac{5}{2}-\frac{ 8^\frac{1}{4}}{2}}{3-\frac{
8^\frac{1}{4}}{2}}{u(\xi_3)} \geq {u(\xi_4)}\geq \frac{\frac{
5}{2}}{2}{u(\xi_3)} .$$

{\bf Step 5.} Let $\xi_5=\exp\left(-\sigma
\sum_{j=1}^nY_j\right)\circ\xi_4=(0_n,0_n,0)\in H_{\xi_4}$. Note
that
$$N(\xi_4)\leq \frac{ 1}{2}R;\ N(\xi_5)=0;\ d_H(\xi_4,\xi_5)=\sigma\sqrt{n}\leq \frac{R}{2}.$$
We may apply Lemma \ref{two-sided-Harnack-lemma} with $c_1=\frac{
1}{2},$ $c_2=0$, $c_3=\frac{1}{2}$ and $c=3$, obtaining
$$\frac{2}{\frac{5}{2}}{u(\xi_4)}\geq {u(\xi_5)}={u(0)} \geq {\frac{3}{\frac{5}{2}}}{u(\xi_4)}.$$
By the Steps 1-5 we conclude that $u(\xi)=0$ if and only if
$u(0)=0$. Therefore, if $u(0)=0$, the arbitrariness of $\xi\in
B_H(0,R)$ shows that $u\equiv 0$ in $B_H(0,R)$.

If $u(0)\neq 0$ then $u<0$ in $B_H(0,R)$, and by multiplying the
estimates from the above five steps, we have that
$$\frac{1}{2}\cdot \frac{3}{4}\cdot \frac{\frac{5}{2}-\frac{ 17^\frac{1}{4}}{2}}{3-\frac{
17^\frac{1}{4}}{2}}\cdot \frac{\frac{5}{2}-\frac{
8^\frac{1}{4}}{2}}{3-\frac{
8^\frac{1}{4}}{2}}\cdot\frac{4}{5}{u(\xi)}\geq {u(0)}\geq 2\cdot
\frac{3-\frac{ 17^\frac{1}{4}}{2}}{\frac{5}{2}-\frac{
17^\frac{1}{4}}{2}}\cdot \frac{3-\frac{
8^\frac{1}{4}}{2}}{\frac{5}{2}-\frac{ 8^\frac{1}{4}}{2}}\cdot
\frac{5}{4}\cdot \frac{6}{5}{u(\xi)}.$$

Repeating the above argument for another point $\zeta\in B_H(0,R)$
and combining the two estimates, it yields that
$$\tilde c^{-1}{u(\xi)}\geq {u(\zeta)}\geq \tilde c{u(\xi)},$$
where $$\tilde c=10\cdot\left(\frac{3-\frac{
17^\frac{1}{4}}{2}}{\frac{5}{2}-\frac{ 17^\frac{1}{4}}{2}}\cdot
\frac{3-\frac{ 8^\frac{1}{4}}{2}}{\frac{5}{2}-\frac{
8^\frac{1}{4}}{2}}\right)^2\approx 30.26,$$ which concludes the
proof.\hfill $\square$ \\

{\it Proof of Proposition \ref{prop-non-positive-H-convex} (second
part)}. Let $\xi_0\in \Omega$ be such that $u(\xi_0)<0$ and fix
$\xi\in \Omega$ arbitrarily. Let
$L=\{(1-\lambda)\xi_0+\lambda\xi:\lambda\in [0,1]\}$ be the
Euclidean segment connecting these two points. From the convexity of
$\Omega$ we conclude that the Euclidean tubular neighborhood around
$L$ with radius $0<r<\min \{{\rm dist}(\xi_0,\partial \Omega),{\rm
dist}(\xi,\partial \Omega)\}$, i.e., $N_L(r)=\{\xi\in \mathbb
H^n:{\rm dist}(\xi,L)<r\},$ is contained in $\Omega.$ [Here, 'dist'
is the Euclidean distance.] Now, we consider the covering
$\bigcup_{\zeta\in L} B_H(\zeta,R_\zeta)$ of the  set $L$ where
$R_\zeta>0$ is such that $B_H(\zeta,3R_\zeta)\subset N_L(r)$ for
every $\zeta\in L.$ By the compactness of $L$, there exists $k\in
\mathbb N$ such that $L\subset \bigcup_{i=1}^k
B_H(\zeta_{\lambda_i},R_{\zeta_{\lambda_i}})$ where
$\zeta_{\lambda_i}=(1-\lambda_i)\xi_0+\lambda_i\xi$ with $0\leq
\lambda_1<...<\lambda_k\leq 1.$ If $k=1$, we are done by
(\ref{harnack-new}), obtaining that $0>\frac{1}{31}u(\xi_0)\geq
u(\xi).$ If $k\geq 2$, since
$B_H(\zeta_{\lambda_i},R_{\zeta_{\lambda_i}})\cap
B_H(\zeta_{\lambda_{i+1}},R_{\zeta_{\lambda_{i+1}}})\neq \emptyset$
for every $i=1,...,k-1,$ we may repeatedly apply (\ref{harnack-new})
on the balls $B_H(\zeta_{\lambda_i},R_{\zeta_{\lambda_i}})$, by
obtaining that $0>\frac{1}{31^k}u(\xi_0)\geq u(\xi).$ \hfill $\square$ \\

\vspace{0.5cm} \noindent {\normalsize\sc  Mathematisches Institute,
Universit\"at Bern,
             Sidlerstrasse 5,
 3012 Bern, Switzerland.
} Email: {\textsf{zoltan.balogh@math.unibe.ch}\\

\noindent {{\normalsize \sc Dipartimento di Matematica e
Applicazioni,
Universit\'a di Milano Bicocca, Via Cozzi 53, 20125 Milano, Italy.}} Email: {\textsf{andrea.calogero@unimib.it}\\

 \noindent {\normalsize \sc Department of Economics, Babe\c s-Bolyai
University,  400591 Cluj-Napoca, Romania.} Email:
{\textsf{alexandru.kristaly@econ.ubbcluj.ro}


\begin{thebibliography}{99}


\bibitem{Aubin-Cellina} J.-P. Aubin, A. Cellina, Differential inclusions. Set-valued maps and viability theory. Springer-Verlag, Berlin, 1984.

\bibitem{Aubin-Fran} J.-P. Aubin, H. Frankowska, Set-Valued Analysis. Systems \& Control: Foundations \& Applications, 2. Birkh\"auser Boston, Inc., Boston, MA, 1990.

\bibitem{forth} Z. M.  Balogh, A. Calogero, A. Krist\'aly, E.
Vecchi, Aleksandrov-type estimates on Carnot groups, in
preparation.

\bibitem{Balogh-Rickly} Z. M.  Balogh, M. Rickly,
Regularity of convex functions on Heisenberg groups. {\it Ann. Sc.
Norm. Super. Pisa Cl. Sci.} (5) 2 (2003), no. 4, 847--868.



\bibitem{Bardi-Dragoni} M. Bardi,  F. Dragoni,
Convexity and semiconvexity along vector fields. {\it Calc. Var.
Partial Differential Equations}, 42 (2011), 405--427.

\bibitem{BLU} A. Bonfiglioli, E. Lanconelli, F. Uguzzoni, Stratified Lie groups and potential theory for their sub-Laplacians. {\it Springer Monographs in Mathematics}, Springer-Verlag, Berlin, 2007.

\bibitem{Caf-1} L. A. Caffarelli, Interior $W^{2,p}-$estimates for solutions of the
Monge-Amp\`ere equation. {\it Ann. Math.} 131 (1990), 135--150.

\bibitem{Caf-2} L. A. Caffarelli, Some regularity properties of solutions of
Monge-Amp\`ere equation. {\it Comm. Pure Appl. Math.} 44 (1991),
965--969.

\bibitem{CaPi2011}
A.~Calogero, R.~Pini,
\newblock Horizontal normal map on the {H}eisenberg group.
\newblock {\em J. Nonlinear Convex Anal.} 12 (2011), no. 2, 287--307.

\bibitem{CaPi2012}
A.~Calogero, R.~Pini,
\newblock {\it c} horizontal convexity on {C}arnot groups.
\newblock {\em J. Convex Analysis} 19 (2012), no. 3, 541--567.

\bibitem{Capogna-Maldonado} L. Capogna, D. Maldonado,  A note on the engulfing property
and the $\Gamma^{1+\alpha}-$regularity of convex functions in Carnot
groups. {\it Proc. Amer. Math. Soc.} 134 (2006), no. 11, 3191--3199.


\bibitem{Cellina} A. Cellina,  Approximations of set-valued functions and fixed point theorems. {\it Ann. Mat. Pura
Appl.} 82 (1969), 17--24.

\bibitem{Cygan} J. Cygan, Subadditivity of homogeneous norms on certain nilpotent Lie groups. {\it Proc. Amer. Math. Soc.} 83, (1981), 69--70.


\bibitem{DaGaNh2003}
D.~Danielli, N.~Garofalo,  D. M. Nhieu,
\newblock Notions of convexity in {C}arnot groups.
\newblock {\em Comm. Anal. Geom.} 11 (2003), no. 2., 263--341.

\bibitem{DaGaNhTo2004}
D.~Danielli, N.~Garofalo,  D. M. Nhieu, F.~Tournier,
\newblock The theorem of Busemann-Feller-Alexandrov in Carnot groups
\newblock{\em Comm. Anal. Geom.} 12 (2004), no. 4., 853--886.


\bibitem{GaTo2005} N.~Garofalo, F.~Tournier,
\newblock New properties of convex functions in the Heisenberg
group.
\newblock {\em Trans. Amer. Math. Soc.} 358 (2006), no. 5, 2011--2055.



\bibitem{Gu2001}
C.~E. Guti\'errez,
\newblock {The Monge-Amp\`ere {E}quation}.
\newblock Birkh\"auser, Boston, MA, 2001.


\bibitem{GuMo-CommPDE} C.~E. Guti\'errez, A.~Montanari,  Maximum and comparison principles for convex functions on the Heisenberg group.
{\it Comm. Partial Differential Equations} 29 (2004), no. 9-10,
1305--1334.

\bibitem{GuMo-Pisa} C.~E. Guti\'errez, A.~Montanari, On the second order derivatives of convex functions on the
   Heisenberg group. {\it Ann. Sc. Norm. Super. Pisa Cl. Sci.} 3 (2004), 349--366.


\bibitem{Hu-Papa} S. Hu, N. S.  Papageorgiou,
Generalizations of Browder's degree theory. {\em Trans. Amer. Math.
Soc.} 347 (1995), no. 1, 233--259.

\bibitem{JuGuMaSt2007} P. Juutinen, G.  Lu, J. J. Manfredi, B. Stroffolini,  {Convex functions on {C}arnot groups}. {\it Rev. Mat.
Iberoam.}
 23  (2007), 191--200.


\bibitem{Lloyd1978} N. G. Lloyd, Degree Theory. Cambridge University Press, 1978.

 \bibitem{GuMaSt2004} G. Lu,  J. J. Manfredi, B.  Stroffolini,  {Convex functions on the Heisenberg group.}
 {\it Calc. Var. Partial Differential Equations},  19 (2004), no. 1,
 1--22.

\bibitem{Magnani} V. Magnani,
Lipschitz continuity, Aleksandrov theorem and characterizations for
$H-$convex functions. {\it Math. Ann.} 334 (2006), 199--233.

\bibitem{Ma-Sci1} V. Magnani, M. Scienza,
Regularity estimates for convex functions in Carnot-Carath\'eodory
spaces. {\it Preprint} (2012).


\bibitem{Ri2006}  M. Rickly, First order regularity of convex functions on {C}arnot groups. {\em J. Geom. Anal. }   16 (2006), no. 4,
679--702.

\bibitem{Ro1969}
R.T. Rockafellar,
\newblock {Convex {A}nalysis}.
\newblock Princeton University Press, 1969.


\bibitem{TrudWang1} N. S. Trudinger, X.-J. Wang,  Hessian measures I. {\it Topol. Methods Nonlinear Analysis,} 10 (1997), no. 2,
225--239.

\bibitem{TrudWang2} N. S. Trudinger, X.-J. Wang, Hessian measures II. {\it Ann. of Math.} (2) 150 (1999), no. 2, 579--604.

\bibitem{TrudWang3} N. S. Trudinger, X.-J. Wang, Hessian measures III. {\it J. Funct. Anal.} 193 (2002), no. 1, 1--23.

\bibitem{TrudZhang} N. S. Trudinger, W. Zhang,  Hessian measures on the
Heisenberg group. {\it J. Funct. Anal.} 264 (2013), no. 10,
2335--2355.

\end{thebibliography}
\end{document}